\documentclass[11pt]{article}

\usepackage{graphicx}

\usepackage{amsthm,amsmath,amssymb}
\usepackage[all,2cell]{xy}
\UseAllTwocells
\usepackage[english]{babel}
\usepackage{subfigure}

\newcommand{\N}{\mathbb{N}}
\newcommand{\Z}{\mathbb{Z}}

\title{About the globular homology of higher dimensional automata}
\author{Philippe Gaucher}

\date{}

\newtheorem{thm}{Theorem}[section]
\newtheorem{prop}[thm]{Proposition}
\newtheorem{propdef}[thm]{Proposition and definition}

\newtheorem{conj}[thm]{Conjecture}
\newtheorem{question}[thm]{Question}
\newtheorem{cor}[thm]{Corollary}
\newtheorem{rem}[thm]{Remark}

\newtheorem{defn}[thm]{Definition}

\newcommand{\be}{\begin{equation}}
\newcommand{\ee}{\end{equation}}
\newcommand{\bea}{\begin{eqnarray}}
\newcommand{\eea}{\end{eqnarray}}
\newcommand{\beas}{\begin{eqnarray*}}
\newcommand{\eeas}{\end{eqnarray*}}

\newcommand{\C}{\mathcal{C}}
\newcommand{\D}{\mathcal{D}}

\newcommand{\bd}{\begin{defn}}
\newcommand{\ed}{\end{defn}}
\newcommand{\bcd}{\begin{defn}}
\newcommand{\ecd}{\end{defn}}
\newcommand{\bp}{\begin{prop}}
\newcommand{\ep}{\end{prop}}
\newcommand{\bth}{\begin{thm}}
\renewcommand{\eth}{\end{thm}}
\newcommand{\bi}{\begin{enumerate}}
\newcommand{\ei}{\end{enumerate}}
\newcommand{\br}{\begin{rem}}
\newcommand{\er}{\end{rem}}
\newcommand{\bpf}{\begin{proof}}
\newcommand{\epf}{\end{proof}}

\newcommand{\p}\times

\newcommand{\iso}{\cong}

\newcommand{\ot}{\otimes}
\newcommand{\de}{\partial}

\renewcommand{\P}{\mathbb{P}}

\renewcommand{\geq}{\geqslant}
\renewcommand{\leq}{\leqslant}
\newcommand{\Cube}{{\underline{cub}}}

\newcommand{\fl}[1]{\ar@{->}[l]_{#1}}
\newcommand{\fr}[1]{\ar@{->}[r]^{#1}}
\newcommand{\fd}[1]{\ar@{->}[d]_{#1}}
\newcommand{\fu}[1]{\ar@{->}[u]^{#1}}
\newcommand{\f}[2]{\ar@{->}[#1]|{#2}}
\newcommand{\ff}[2]{\ar@2{->}[#1]|{#2}}

\newcommand{\HR}{H\!R}
\newcommand{\CR}{C\!R}
\newcommand{\HF}{H\!F}
\newcommand{\CF}{C\!F}
\newcommand{\ev}{e\!v}

\newcommand{\comp}{C\!omp(Ab)}
\newcommand{\id}{I\!d}

\begin{document}

\maketitle

\hbox{ }
\vspace{-1cm}

\begin{center}\textbf{R\'esum\'e}  \end{center}
\begin{list}{}
{\leftmargin=1.5cm}
  \item On introduit un nouveau nerf simplicial d'automate parall\`ele dont
  l'homologie simpliciale d\'ecal\'ee de un fournit une nouvelle
  d\'e\-fi\-ni\-tion de l'homologie globulaire. Avec cette nouvelle
  d\'e\-fi\-ni\-tion, les inconv\'enients de la construction de \cite{Gau}
  disparaissent.  De plus les importants morphismes qui associent \`a
  tout globe les zones correspondantes de branchements et de
  confluences de chemins d'e\-x\'e\-cu\-tion deviennent ici des morphismes
  d'ensembles simpliciaux.
\end{list}

\tableofcontents

\section{Introduction}\label{intro}

One of the contributions of \cite{HDA} is the introduction of two
homology theories as a starting point for studying branchings and
mergings in higher dimensional automata (HDA) from an homological
point of view. However these homology theories had an important
drawback : roughly speaking, they were not invariant by subdivisions
of the observation. Later in \cite{Gau}, using a model of concurrency
by strict globular $\omega$-categories borrowed from \cite{Pratt}, two
new homology theories are introduced : the negative and positive
corner homology theories $H^-$ and $H^+$, also called the branching
and the merging homologies.  It is proved in \cite{Coin} that they
overcome the drawback of Goubault's homology theories.

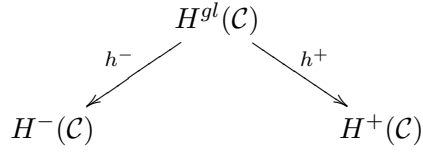
\begin{figure}
\[\xymatrix{& {H^{gl}(\C)}\ar@{->}[ld]_{h^-}\ar@{->}[rd]^{h^+}&\\
{H^{-}(\C)}&&{H^{+}(\C)}}\]
\caption{Associating to any globe its two corners}
\label{fundamental}
\end{figure}

Another idea of \cite{Gau} is the construction of a diagram of abelian
groups like in Figure~\ref{fundamental}, where $H_*^{gl}$ is a new
homology theory called the globular homology.

Geometrically, the non-trivial cycles of the globular homology must
correspond to the oriented empty globes of $\C$, and the non-trivial
cycles of the branching (resp. the merging) homology theory must
correspond to the branching (resp. merging) areas of execution paths.
And the morphisms $h^-$ and $h^+$ must associate to any globe its
corresponding branching area and merging area of execution paths.
Many potential applications in computer science of these morphisms are
put forward in
\cite{Gau}.

Globular homology was therefore created in order to fulfill two
conditions :

\begin{itemize}
\item Globular homology must take place in a diagram of abelian
  groups like in Figure~\ref{fundamental}. And the geometric meaning
  of $h^-$ and $h^+$ must be exactly as above described.
\item Globular homology must be an invariant of HDA with respect to
  reasonable deformations of HDA, that is of the corresponding
  $\omega$-category.
\end{itemize}

What is a reasonable deformation of HDA was not yet very clear in
\cite{Gau}.  This question is discussed with much more details in
\cite{ConcuToAlgTopo}.

The \textit{old globular homology} (i.e. the construction exposed
in \cite{Gau}) satisfied the first condition, and the second one
was supposed to be satisfied by definition (cf. Definition 8.2 of
two homotopic $\omega$-categories in \cite{Gau}), even if some
problems were already mentioned, particularly the non-vanishing
of the ``old'' globular homology of $I^3$, and more generally of
$I^n$ for any $n\geq 1$ in strictly positive dimension.

This latter problem is disturbing because the $n$-cube $I^n$ (i.e.
the corresponding automaton which consists of $n$ $1$-transitions
carried out at the same time) can be deformed by crushing all the
$p$-faces with $p>1$ into an $\omega$-category which has only
$0$-morphisms and $1$-morphisms and because the globular homology is
supposed to be an invariant by such deformations.  The philosophy
exposed in \cite{ConcuToAlgTopo} tells us similar things : using
S-deformations and T-deformations, the $n$-cube and the oriented line
must be the same up to homotopy, and therefore must have the same
globular homology.

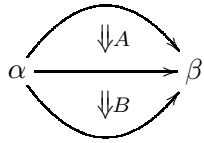
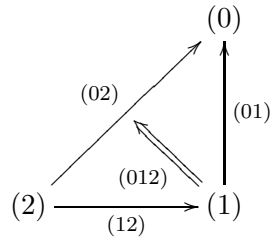
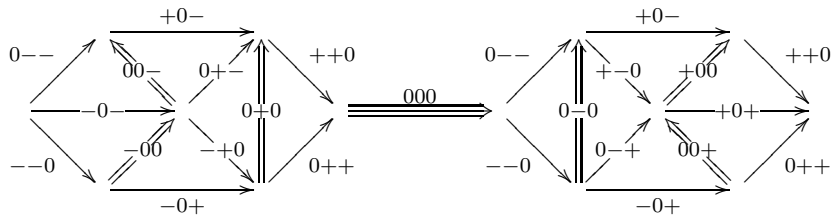
\begin{figure}
\begin{center}
\subfigure[Composition of two $2$-morphims]{\label{composition-of-2trans}
\xymatrix{&&\\
{\alpha} \rruppertwocell<10>{A} \rrlowertwocell<-10>{B} \ar[rr]&& {\beta}
}}
\hspace{2cm}
\subfigure[The  $\omega$-category $\Delta^2$]{\label{2simpl}
\xymatrix{&& {(0)}\\&&\\{(2)}\ar@{->}[rruu]^{(02)}\ar@{->}[rr]_{(12)}&&{(1)}\ar@{->}[uu]_{(01)}
\ar@2{->}[lu]^{(012)}}}
\subfigure[The $\omega$-category $I^3$]{\label{I3}
\xymatrix{
&\ar@{->}[rr]^{+0-}& &\ar@{->}[rd]^{++0}&&&
&\ar@{->}[dr]|{+-0}\ar@{->}[rr]^{+0-}&&\ar@{->}[dr]^{++0}&\\
\ar@{->}[ru]^{0--}\ar@{->}[rr]|{-0-}\ar@{->}[rd]_{--0}&&\ar@{->}[ru]|{0+-}\ar@{->}[dr]|{-+0}\ff{lu}{00-}&&
\ar@3{->}[rr]^{000}&&
\ar@{->}[ur]^{0--}\ar@{->}[dr]_{--0}&&\ff{ur}{+00}\ar@{->}[rr]|{+0+}&&\\
&\ar@{->}[rr]_{-0+}\ff{ru}{-00}&&\ar@{->}[ru]_{0++}\ff{uu}{0+0}&&&
&\ar@{->}[ru]|{0-+}\ar@{->}[rr]_{-0+}\ff{uu}{0-0}&&\ff{ul}{00+}\ar@{->}[ru]_{0++}\\
}}
\end{center}
\caption{Some $\omega$-categories (a $k$-fold arrow symbolizes a k-morphism)}
\end{figure}

The non-vanishing of the second globular homology group of $I^3$ (see
Figure~\ref{I3}) is due for instance to the $2$-dimensional globular
cycle

\beas
&&\left(R(-00) *_0 R(0++)\right)*_1 \left(R(-0-) *_0
R(0+0)\right) \\ &&-
\left(R(-00) *_0 R(0++)\right) - \left(R(-0-) *_0 R(0+0)\right)
\eeas

It is the reason why it was suggested in \cite{Gau}
to add the relation $A*_1 B=A+B$ at least to the $2$-dimensional stage
of the old globular complex.

But there is then no reason not to add the same relation in the rest of the
definition of the old globular complex. For example, if we take the quotient
of the old globular complex by the relation $A*_1 B=A+B$ for any pair
$(A,B)$ of $2$-morphisms, then the
$\omega$-category defined as the free $\omega$-category generated by
the globular set generated by two $3$-morphisms $A$ and $B$ such that
$t_1 A=s_1 B$ gives rise to a $3$-dimensional globular cycle $A
*_1 B -A-B$ because $s_2(A *_1 B -A-B)=s_2 A *_1 s_2 B -s_2 A -s_2
B=0$ and $t_2(A *_1 B-A-B)=t_2 A *_1 t_2 B -t_2 A -t_2 B=0$.  So putting
the relation $A*_1 B -A-B=0$ in the old globular complex for any pair
of morphisms $(A,B)$ of the same dimension sounds necessary.
Similar considerations starting from the calculation of the $(n-1)$-th
globular homology group of $I^n$ entail the relations $A *_n B-A-B$
for any $n\geq 1$ and for any pair $(A,B)$ of $p$-morphisms with $p
\geq n+1$ in the old globular chain complex.

The \textit{formal globular homology} of
Definition~\ref{new_glob} is exactly equal to the quotient of the
old globular complex by these missing relations. It is
conjectured (see conjecture~\ref{thin_glob}) that this homology
theory will coincide for free $\omega$-categories generated by
semi-cubical sets with the homology theory of
Definition~\ref{def_correcte}, this latter being the simplicial
homology of the globular simplicial nerve $\mathcal{N}^{gl}$
shifted by one.

We claim that Definition~\ref{def_glob_nerve} (and its simplicial
homology shifted by one) cancels the drawback of the old globular
homology at least for the following reasons :

\begin{itemize}
\item It is noticed in \cite{Gau} that both corner homologies come
  from the simplicial homology of two augmented simplicial nerves
  $\mathcal{N}^-$ and $\mathcal{N}^+$ ; there exists one and only
  natural transformation $h^-$ (resp.  $h^+$) from $\mathcal{N}^{gl}$
  to $\mathcal{N}^-$ (resp. $\mathcal{N}^+$) preserving the interior
  labeling (Theorem~\ref{hmoins}).
\item In homology, $h^-$ and $h^+$ induce two natural linear maps from
  $H_*^{gl}$ to resp. $H_*^-$ and $H_*^+$ which do exactly what we want.
\item The globular homology (formal or not) of $I^n$ vanishes in
  strictly positive dimension for any $n\geq 0$. The globular homology
  of $\Delta^n$ (the $n$-simplex) and of $2_n$ (the free
  $\omega$-category generated by one $n$-di\-men\-sio\-nal morphism) as
  well.
\item Using Theorem~\ref{relation-old_new} explaining the exact mathematical
link between the old construction and the new one, one sees that one
does not lose the possible applications in computer science pointed out
in \cite{Gau}.
\item The new globular homology, as well as the new globular cut are
  invariant by S-deformations, that is intuitively by contraction and
  dilatation of homotopies between execution paths. We will see
  however that it is not invariant by T-deformations, that is by
  subdivision of the time, as the old definition and
this problem will be a little bit discussed.
\end{itemize}

This paper is two-fold. The first part introduces the new
material. The second part justifies the new definition of the
globular homology.

After Section~\ref{convention} which recalls some conventions
and some elementary facts about strict globular $\omega$-categories
(non-contracting or not) and about simplicial sets, the setting of
\textit{simplicial cuts} of non-contracting $\omega$-categories and that of
\textit{regular cuts} are introduced. The first notion allows to enclose
the new globular nerve of this paper and both corner nerves in one
unique formalism. The notion of regular cuts gives an axiomatic
framework for the generalization of the notion of negative and
positive folding operators of
\cite{Coin}. Section~\ref{example_corner} is an illustration of the previous
new notions on the case of corner nerves. In the same section, some
non-trivial facts about negative folding operators are recalled.
Section~\ref{glob_cut} provides the definition of the globular
nerve of a non-contracting $\omega$-category.

The organization of the rest of the paper follows the preceding
explanations. First in Section~\ref{morphisms}, the morphisms $h^-$
and $h^+$ are constructed. Section~\ref{regularity} proves that the
globular cut is regular. In particular, we get the globular folding
operators. Section~\ref{vanish_In} proves the vanishing of the
globular homology of the $n$-cube, the $n$-simplex and the free $\omega$-category
generated by one $n$-morphism. At last Section~\ref{relation} makes
explicit the exact relation between the new globular homology and the
old one. Section~\ref{deformation} speculates about deformations of
$\omega$-categories considered as a model of HDA and the construction
of the bisimplicial set of \cite{ConcuToAlgTopo} is detailed.

\section{Conventions and notations}\label{convention}

\subsection{Globular $\omega$-category and cubical set}

For us, an
$\omega$-category will be a strict globular $\omega$-category with
morphisms of finite dimension. More precisely (see
\cite{Brown-Higgins0} \cite{oriental} \cite{Tensor_product} for more
details) :

\bd A \textit{$1$-category}
is a pair $(A,(*,s,t))$ satisfying the following
axioms\thinspace:
\begin{enumerate}
\item $A$ is a set
\item $s$ and $t$ are set maps from $A$ to $A$ respectively called the
source map and the target map
\item for $x,y\in A$, $x*y$ is defined as soon as $tx=sy$
\item $x*(y*z)=(x*y)*z$ as soon as both members of the equality exist
\item $sx *x=x*tx=x$, $s(x*y)=sx$ and $t(x*y)=ty$ (this implies
$ssx=sx$ and $ttx=tx$).
\end{enumerate}
\ed

\bd\label{2categories}
A \textit{$2$-category}
is a triple $(A,(*_0,s_0,t_0),(*_1,s_1,t_1))$ such that
\begin{enumerate}
\item both pairs $(A,(*_0,s_0,t_0))$ and $(A,(*_1,s_1,t_1))$ are $1$-categories
\item $s_0 s_1=s_0 t_1=s_0$, $t_0 s_1=t_0 t_1=t_0$, and for $i\geq j$,
$s_i s_j=t_i s_j=s_j$ and $s_i t_j=t_i t_j=t_j$ (Globular axioms)
\item $(x*_0y)*_1 (z*_0 t)=(x*_1z)*_0 (y*_1 t)$ (Godement axiom or interchange law)
\item if $i\neq j$, then $s_i (x*_jy) =s_i x *_j s_i y$ and
$t_i (x*_jy)= t_i x *_j t_i y$.
\end{enumerate}
\ed

\bd\label{omega_categories}  A \textit{globular $\omega$-category} $\C$
is a set $A$ together with a
family $(*_n,s_n,t_n)_{n\geq 0}$ such that
\begin{enumerate}
\item for any $n\geq 0$, $(A,(*_n,s_n,t_n))$ is a $1$-category
\item for any $m,n\geq 0$ with $m<n$,
$(A,(*_m,s_m,t_m),(*_n,s_n,t_n))$ is a $2$-category
\item for any $x\in A$, there exists $n\geq 0$ such that $s_nx=t_nx=x$
(the smallest of these $n$ is called the dimension of $x$).
\end{enumerate}
\ed

A $n$-dimensional element of $\C$  is called a $n$-morphism.
A $0$-morphism is also called a \textit{state} of $\C$, and a $1$-morphism
an \textit{arrow}. If $x$ is a morphism of an $\omega$-category $\C$, we call
$s_n(x)=d_n^-(x)$ the \textit{$n$-source} of $x$ and $t_n(x)=d_n^+(x)$ the
\textit{$n$-target} of $x$.  The category of all $\omega$-categories (with the
obvious morphisms) is denoted by $\omega Cat$. The corresponding
morphisms are called $\omega$-functors. The set of
morphisms of $\C$ of dimension at most $n$ is denoted by $tr^n \C$ ;
the set of morphisms of $\C$ of dimension exactly $n$ is denoted by
$\C_n$.

Sometime we will use the terminology \textit{initial state}
(resp. \textit{final state}) for a state $\alpha$ which is not the
$0$-target (resp. the $0$-source) of a $1$-morphism.

\bd\label{def_cubique}\cite{Brown_cube} \cite{cube}
A \textit{cubical set} consists of
\begin{itemize}
\item a family of
sets $(K_n)_{n\geq 0}$
\item a family of face maps
$\xymatrix@1{{K_n}\fr{\de_i^\alpha} &{K_{n-1}}}$ for $\alpha\in\{-,+\}$
\item a family of degeneracy maps
$\xymatrix@1{K_{n-1}\fr{\epsilon_i}&K_{n}}$ with $1\leq
i\leq n$
\end{itemize}
which satisfy the following relations
\begin{enumerate}
\item $\de_i^\alpha \de_j^\beta = \de_{j-1}^\beta \de_i^\alpha$
for all $i<j\leq n$ and $\alpha,\beta\in\{-,+\}$
\item $\epsilon_i\epsilon_j=\epsilon_{j+1}\epsilon_i$
for all $i\leq j\leq n$
\item $\de_i^\alpha \epsilon_j=\epsilon_{j-1}\de_i^\alpha$       for
  $i<j\leq n$
and $\alpha\in\{-,+\}$
\item $\de_i^\alpha \epsilon_j=\epsilon_{j}\de_{i-1}^\alpha$   for
  $i>j\leq n$ and
$\alpha\in\{-,+\}$
\item $\de_i^\alpha \epsilon_i=Id$
\end{enumerate}
A family $(K_n)_{n\geq 0}$ only equipped with a family of face maps
$\de_i^\alpha$ satisfying the same axiom as above is called a \textit{semi-cubical set}.
\ed

\bd
The corresponding category of cubical sets, with an obvious
definition of its morphisms, is isomorphic to the category of
presheaves $Sets^{\square^{op}}$ over a small category $\square$. The
corresponding category of semi-cubical sets , with an obvious
definition of its morphisms, is isomorphic to the category of
presheaves $Sets^{{\square^{semi}}^{op}}$ over a small category $\square^{semi}$.
\ed

In a simplicial set, the face maps are always denoted by $\de_i$, the
degeneracy maps by $\epsilon_i$. Here are the other conventions about
simplicial sets (see for example \cite{May} for further information) :

{\small
\begin{enumerate}
\item $Sets$ : category of sets
\item $Sets^{\Delta^{op}}$ : category of simplicial sets
\item $\comp$ : category of chain complexes of abelian groups
\item $C(A)$ : unnormalized chain complex of the simplicial set $A$
\item $H_*(A)$ : simplicial homology of a simplicial set $A$
\item $Ab$ : category of abelian groups
\item $\id$  : identity map
\item $\Z S$ : free abelian group generated by the set $S$
\end{enumerate}}

HDA means \textit{higher dimensional automaton}. In this paper, this
is another term for \textit{semi-cubical set}, or the corresponding free
$\omega$-category generated by it.

Various homology theories (see the diagram of
Theorem~\ref{relation-old_new}) will appear in this paper.  It is
helpful for the reader to keep in mind that the total homology of a
semi-cubical set is used nowhere in this work.

\subsection{Non-contracting $\omega$-category}

Let $\C$ be an $\omega$-category. We want to define an
$\omega$-category $\P\C$ ($\P$ for path) obtained from $\C$ by
removing the $0$-morphisms, by considering the $1$-mor\-phisms of $\C$
as the $0$-morphisms of $\P\C$, the $2$-morphisms of $\C$ as the
$1$-morphisms of $\P\C$ etc. with an obvious definition of the source
and target maps and of the composition laws (this new
$\omega$-category is denoted by $\C[1]$ in \cite{ConcuToAlgTopo}). The
map $\P:\C\mapsto \P\C$ does not induce a functor from $\omega Cat$ to
itself because $\omega$-functors can contract $1$-morphisms and
because with our conventions, a $1$-source or a $1$-target can be
$0$-dimensional. Hence the following definition

\begin{propdef}\label{noncontract}
For a globular $\omega$-category $\C$, the following
assertions are equivalent :
\begin{enumerate}
\item[(i)] $\P\C$ is an $\omega$-category ; in other terms, $*_i$, $s_i$
and $t_i$ for any $i\geq 1$ are internal to $\P\C$ and we can set
$*_i^{\P\C}=*_{i+1}^{\C}$, $*_i^{\P\C}=*_{i+1}^{\C}$ and $*_i^{\P\C}=*_{i+1}^{\C}$ for any
$i\geq 0$.
\item[(ii)] The maps $s_1$ and $t_1$ are
non-contracting,
that is if $x$ is of strictly positive dimension, then $s_1x$ and $t_1x$ are
$1$-dimensional (a priori, one can only say that $s_1x$ and $t_1x$ are
of
dimension lower or equal than $1$)
\end{enumerate}
If Condition (ii) is satisfied, then one says that $s_1$ and
$t_1$
are non-contracting and that $\C$ is \textit{non-contracting}.
\end{propdef}

\bpf Suppose $s_1$ and $t_1$ non-contracting. Let $x$ and $y$ be
two morphisms of strictly positive dimension and $p\geq 1$. Then $s_1
s_p x=s_1x$ therefore $s_px$ cannot be $0$-dimensional. If $x *_p y$
then $s_1(x*_py)=s_1x$ if $p=1$ and if $p>1$ for two different
reasons. Therefore $x*_py$ cannot be $0$-dimensional as soon as $p\geq 1$.
\epf

\bd\label{noncontractant} Let $f$ be an $\omega$-functor from $\C$ to
$\D$.  The morphism $f$ is \textit{non-contracting} if for any
$1$-dimensional $x\in \C$, the morphism $f(x)$ is a $1$-dimensional
morphism of $\D$ (a priori, $f(x)$ could be either $0$-dimensional
or $1$-dimensional).  \ed

\bd The category of non-contracting $\omega$-categories with the
non-contracting $\omega$-functors is denoted by $\omega Cat_1$.  \ed

Notice that in \cite{Gau}, the word ``non-$1$-contracting'' is used
instead of simply ``non-contrac\-ting''. Since \cite{ConcuToAlgTopo},
the philosophy behind the idea of deforming the $\omega$-categories
viewed as models of HDA is better understood. In particular, the idea
of not contracting the morphisms is relevant only for $1$-dimensional
morphisms. So the ``1'' in ``non-$1$-contracting'' is not anymore
necessary.

\bd Let $\C$ be a non-contracting $\omega$-category. Then the $\omega$-category
$\P\C$ above defined is called the \textit{path $\omega$-category} of $\C$.
The map $\C \mapsto \P\C$ induces a functor from $\omega Cat_1$ to
$\omega Cat$. \ed

Here is a fundamental example of non-contracting $\omega$-category.
Consider a semi-cubical set $K$ and consider the free $\omega$-category
$\Pi(K):=\int^{\underline{n}\in \square} K_n .I^n$ generated by it
where
\begin{itemize}
\item $I^n$ is the free $\omega$-category generated by the faces of
  the $n$-cube, whose construction is recalled in
  Section~\ref{example_corner}.
\item the integral sign denotes the coend construction and
$K_n.I^n$ means the sum of ``cardinal of $K_n$''
copies of $I^n$ (cf. \cite{cat} for instance).
\end{itemize}
Then one has

\bp For any semi-cubical set $K$, $\Pi(K)$ is a non-con\-trac\-ting
$\omega$-category.  The functor
$\Pi:Sets^{{\square^{semi}}^{op}}\rightarrow \omega Cat $ from the
category of semi-cubical sets to that of $\omega$-categories yields a
functor from $Sets^{{\square^{semi}}^{op}}$ to the category of
non-contracting $\omega$-categories $\omega Cat_1$. \ep

\bpf The  characterization of Proposition~\ref{noncontract} gives
the solution. \epf

\section{Cut of globular higher dimensional categories}\label{cut}

Before introducing the globular nerve of an $\omega$-category, let us
introduce the formalism of \textit{regular simplicial cuts} of
$\omega$-categories. The notion of \textit{simplicial cuts} enables us to put
together in the same framework both corner nerves constructed in
\cite{Gau,Coin} and the new globular nerve of Section~\ref{glob_cut}.  The
notion of \textit{regular  cuts} enables to generalize the notion of
negative (resp. positive) folding operators associated to the
branching (resp. merging) nerve (cf. \cite{Coin}). It is also an
attempt to finding a way of characterizing these three nerves. There are
no much more things known about this problem.

\bd\cite{triple} An augmented simplicial set is a simplicial
set
\[((X_n)_{n\geq 0}, (\de_i:X_{n+1}\longrightarrow X_n)_{0\leq i\leq n+1},
(\epsilon_i:X_{n}\longrightarrow X_{n+1})_{0\leq i\leq n})\] together
with an additional set $X_{-1}$ and an additional map $\de_{-1}$ from
$X_0$ to $X_{-1}$ such that $\de_{-1}\de_0=\de_{-1}\de_1$. A morphism
of augmented simplicial set is a map of $\N$-graded sets which
commutes with all face and degeneracy maps. We denote by
$Sets^{\Delta^{op}}_+$ the category of augmented simplicial sets.  \ed

The ``chain complex'' functor of an augmented simplicial set $X$ is
defined by $C_n(X)=\Z X_n$ for $n\geq -1$ endowed with the simplicial
differential map (denoted by $\de$)
in positive dimension and the map $\de_{-1}$ from
$C_0(X)$ to $C_{-1}(X)$. The ``simplicial homology''
functor $H_*$ from the category of augmented simplicial sets
$Sets^{\Delta^{op}}_+$ to the category of abelian groups $Ab$ is
defined as the usual one for $*\geq 1$ and by setting
$H_0(X)=Ker(\de_{-1})/Im(\de_0-\de_1)$ and
$H_{-1}(X)=\Z X_{-1}/Im(\de_{-1})$ whenever $X$ is an augmented
simplicial set.

\bd A \textit{(simplicial) cut} is a functor $\mathcal{F}:\omega Cat_1 \rightarrow
Sets^{\Delta^{op}}_+$ together with a family $\ev=(\ev_n)_{n\geq 0}$
of natural transformations $\ev_n:F_n\longrightarrow tr^n \P$ where
$F_n$ is the set of $n$-simplexes of $\mathcal{F}$.
A morphism of cuts from $(\mathcal{F},\ev)$ to
$(\mathcal{G},\ev)$ is a natural transformation of
functors $\phi$ from $\mathcal{F}$ to $\mathcal{G}$ which makes the following diagram
commutative for any $n\geq 0$ :
\[\xymatrix{{\mathcal{F}_n} \fd{\phi_n} \fr{\ev_n}& {tr^n \P}\\ {\mathcal{G}_n} \ar@{->}[ru]_{\ev_n}&}\]
\ed

The terminology of ``cuts'' is borrowed from \cite{cyl}.  It will be
explained later : cf. the explanations around Figure~\ref{2simplglob} and
also Section~\ref{deformation}.

There is no
ambiguity to denote all $\ev_n$ by the same notation $\ev$ in the
sequel. The map $\ev$ of $\N$-graded sets is called the
\textit{evaluation map} and a cut $(\mathcal{F},\ev)$ will be always
denoted by $\mathcal{F}$.

If $\mathcal{F}$ is a functor from $\omega Cat_1$ to $Sets^{\Delta^{op}}_+$, let
$C_{n+1}^\mathcal{F}(\C):=C_n(\mathcal{F}(\C))$ and let $H_{n+1}^\mathcal{F}$ be the
corresponding homology theory for $n\geq -1$.

Let $M_n^\mathcal{F}:\omega Cat_1\longrightarrow Ab$ be the functor defined as
follows : the group $M_{n}^\mathcal{F}(\C)$ is the subgroup generated by the
elements $x\in \mathcal{F}_{n-1}(\C)$ such that $\ev(x)\in tr^{n-2}\P\C$ for
$n\geq 2$ and with the convention $M_0^\mathcal{F}(\C)=M_1^\mathcal{F}(\C)=0$ and the
definition of $M_n^\mathcal{F}$ is obvious on non-contracting
$\omega$-functors. The elements of $M_{*}^\mathcal{F}(\C)$ are called \textit{thin}.

Let $\CR^\mathcal{F}_n:\omega Cat_1\longrightarrow \comp$ be the functor
defined by $\CR^\mathcal{F}_n:=C_n^\mathcal{F}/(M_n^\mathcal{F}+\de M_{n+1}^\mathcal{F})$ and endowed with the
differential map $\de$.  This chain complex is called the \textit{reduced
complex} associated to the cut $\mathcal{F}$ and the corresponding homology is
denoted by $\HR^\mathcal{F}_*$ and is called the \textit{reduced homology} associated to
$\mathcal{F}$. A morphism of cuts from $\mathcal{F}$ to $\mathcal{G}$ yields natural morphisms from
$H_*^\mathcal{F}$ to $H_*^\mathcal{G}$ and from $\HR_*^\mathcal{F}$ to $\HR_*^\mathcal{G}$. There is also a
canonical natural transformation $R^\mathcal{F}$ from $H_*^\mathcal{F}$ to $\HR_*^\mathcal{F}$, functorial with
respect to $\mathcal{F}$, that is making the following diagram commutative :
\[
\xymatrix{{H_*^\mathcal{F}}\fr{R^\mathcal{F}}\fd{} & {\HR_*^\mathcal{F}}\fd{}\\ {H_*^\mathcal{G}}\fr{R^\mathcal{G}}& {\HR_*^\mathcal{G}}}
\]

\bd\label{def_regular}
A cut $\mathcal{F}$ is \textit{regular} if and only if it satisfies the following
properties :
\begin{enumerate}
\item\label{regular0} For any $\omega$-category $\C$, the set $\mathcal{F}_{-1}(\C)$ only depends on $tr^0\C=\C_0$ :
i.e. for any $\omega$-categories $\C$ and $\D$, $\C_0=\D_0$ implies $\mathcal{F}_{-1}(\C)=\mathcal{F}_{-1}(\D)$.
\item\label{regular1} $\mathcal{F}_0:=tr^0 \P$.
\item \label{regular1.5} $\ev \circ \epsilon_i=\ev$.
\item\label{regular2} for any natural transformation of functors $\mu $ from $\mathcal{F}_{n-1}$ to
 $\mathcal{F}_{n}$ with $n\geq 1$, and for any natural map $\square$ from
$tr^{n-1}\P$ to $\mathcal{F}_{n-1}$ such that
$\ev\circ \square=\id_{tr^{n-1}\P}$, there exists one and only one
natural transformation $\mu .\square$ from $tr^n\P$ to $\mathcal{F}_{n}$ such
that the following diagram commutes
\[\xymatrix{\ar@/^20pt/[rr]^{\id_{tr^n\P}}{tr^n \P}\fr{\mu .\square}& {\mathcal{F}_n}\fr{\ev_n}&{tr^n \P}\\
\ar@/_20pt/[rr]_{\id_{tr^{n-1}\P}}{tr^{n-1}\P}\fu{i_{n}}\fr{\square}&{\mathcal{F}_{n-1}}\fr{\ev_{n-1}}\fu\mu &{tr^{n-1}\P}\fu{i_{n}}}\]
where $i_{n}$ is the canonical inclusion functor from
$tr^{n-1}\P$ to $tr^n \P$.
\item\label{regular3}  let $\square_1^\mathcal{F}:=\id_{\mathcal{F}_0}$ and
$\square_n^\mathcal{F}:=\epsilon_{n-2}.\dots\epsilon_0 . \square_1^\mathcal{F}$ a natural
transformation from $tr^{n-1}\P$ to $\mathcal{F}_{n-1}$ for
$n\geq 2$ ; then the natural transformations $\de_i\square_n^\mathcal{F}$ for $0\leq i \leq n-1$
from $tr^{n-1}\P$ to $\mathcal{F}_{n-2}$ satisfy the following properties
\begin{enumerate}
\item\label{regular3.01} $\left\{\ev \de_{n-2}\square_n^\mathcal{F}, \ev \de_{n-1}\square_n^\mathcal{F}\right\}=\left\{s_{n-1},t_{n-1}\right\}$.
\item\label{regular3.02} if for some $\omega$-category $\C$ and some $u\in \C_n$, $\ev\de_i\square_n^\mathcal{F}(u)=d_p^\alpha u$
for some $p\leq n$ and for some $\alpha\in\{-,+\}$, then $\de_i\square_n^\mathcal{F}=\de_i\square_n^\mathcal{F} d_p^\alpha $.
\end{enumerate}
\item\label{regular3.1} Let  $\Phi_n^\mathcal{F}:=\square_n^\mathcal{F}\circ \ev$ be a natural
transformation from $\mathcal{F}_{n-1}$ to itself\thinspace; then $\Phi_n^\mathcal{F}$
induces the identity natural transformation  on $\CR_{n}^\mathcal{F}$.
\item\label{regular4} if $x$, $y$ and $z$ are three elements of $\mathcal{F}_n(\C)$,
and if $\ev(x)*_p \ev(y)=\ev(z)$ for some $1\leq p\leq n$, then
$x+y=z$ in $\CR_{n+1}^{\mathcal{F}}(\C)$ and in a functorial way.
\end{enumerate}
If $\mathcal{F}$ is a regular cut, then the natural transformation $\Phi_n^\mathcal{F}$ is
called the $n$-dimensional folding operator of the cut $\mathcal{F}$. By
convention, one sets $\square_0^\mathcal{F}=\id_{\mathcal{F}_{-1}}$ and
$\Phi_0^\mathcal{F}=\id_{\mathcal{F}_{-1}}$. There is no ambiguity to set
$\Phi^\mathcal{F}(x):=\Phi^\mathcal{F}_{n+1}(x)$ for $x\in \mathcal{F}_n(\C)$ for some
$\omega$-category $\C$. So $\Phi^\mathcal{F}$ defines a natural transformation,
and even a morphism of cuts, from $\mathcal{F}$ to itself. However beware of the
fact that there is really an ambiguity in the notation $\square^\mathcal{F}$ :
so this latter will not be used.
\ed

Condition~\ref{regular1.5} tells us that the $\epsilon_i$ operations
are really degeneracy maps.  Condition~\ref{regular2} ensures the
existence and the uniqueness of the folding operator associated to the
cut.

Condition~\ref{regular3} tells us several things. A priori, a natural
transformation like $\ev \de_i\square_n^\mathcal{F}$ from $tr^{n-1}\P$ to
$tr^{n-2}\P$ is necessarily of the form $d_p^\alpha$ for some $p\leq
n-1$ and for some $\alpha\in\{-,+\}$. Indeed consider the free
$\omega$-category $2_n(A)$ generated by some $n$-morphism $A$. Then
$\ev \de_i\square_n^\mathcal{F}(A)\in 2_n(A)$ and therefore $\ev
\de_i\square_n^\mathcal{F}(A) =d_p^\alpha(A)$ for some $p$ and some $\alpha$.
By naturality, this implies that $\ev \de_i\square_n^\mathcal{F}= d_p^\alpha$.
If $0\leq i< n-2$, then
\begin{alignat*}{2}
\ev \de_i \square_n^\mathcal{F} &= \ev \de_i \square_n^\mathcal{F} d_{n-1}^\beta &&\ \hbox{ for some $\beta\in\{-,+\}$}\\
 &= \ev \de_i \square_n^\mathcal{F} i_{n-1} d_{n-1}^\beta &&\\
&=    \ev    \de_i \epsilon_{n-2}\square_{n-1}^\mathcal{F} d_{n-1}^\beta&&\ \hbox{ by construction of $\square_n^\mathcal{F}$}\\
&= \ev\epsilon_{n-3}\de_i \square_{n-1}^\mathcal{F} d_{n-1}^\beta &&\\
&= \ev\de_i \square_{n-1}^\mathcal{F} d_{n-1}^\beta&&\ \hbox{ by rule~\ref{regular1.5}}\\
&= d_p^\alpha d_{n-1}^\beta&&\ \hbox{ for some $p\leq n-2$}\\
&= d_p^\alpha &&
\end{alignat*}
Therefore $\de_i \square_n^\mathcal{F}$ is thin.
Now if $n-2\leq i \leq n-1$, then
\begin{alignat*}{2}
\ev \de_i \square_n^\mathcal{F} &= \ev \de_i \square_n^\mathcal{F} d_{n-1}^\beta &&\ \hbox{ for some $\beta\in\{-,+\}$}\\
 &= \ev \de_i \square_n^\mathcal{F} i_{n-1} d_{n-1}^\beta &&\\
&=    \ev    \de_i \epsilon_{n-2}\square_{n-1}^\mathcal{F} d_{n-1}^\beta&&\ \hbox{ by construction of $\square_n^\mathcal{F}$}\\
&=\ev \square_{n-1}^\mathcal{F} d_{n-1}^\beta &&\\
&= d_{n-1}^\beta &&\ \hbox{ by construction of $\square_n^\mathcal{F}$}
\end{alignat*}
Therefore $\left\{\ev \de_{n-2}\square_n^\mathcal{F}, \ev
  \de_{n-1}\square_n^\mathcal{F}\right\}\subset\left\{s_{n-1},t_{n-1}\right\}$ 
always holds.
Condition~\ref{regular3} states more precisely that these latter sets
are actually equal.  In other terms, the operator
$\square_n^\mathcal{F}$ concentrates the ``weight'' on the faces
$\de_{n-2}\square_n^\mathcal{F}$ and $\de_{n-1}\square_n^\mathcal{F}$.

Condition~\ref{regular3.1} explains the link between the thin
elements of the cut and the folding operators.  Intuitively, the
folding operators move the labeling of the elements of the cuts
in a \textit{canonical position} without changing the total sum
on the source and target sides. What is exactly this
\textit{canonical position} is precisely described by
Proposition~\ref{canonical_position}. Conditions~\ref{regular3}
and \ref{regular4} ensure that by moving the labeling of an
element, we stay in the same equivalence class modulo thin
elements.

Now here are some trivial remarks about regular cuts :
\begin{itemize}
\item Let $f$ be a natural set map from $tr^0 \P\C=\C_1$ to itself. Let $2_1$ be the
$\omega$-category generated by one $1$-morphism $A$. Then necessarily
$f(A)=A$ and therefore $f=\id$. So the above axioms imply that
$\ev_0=\id$.

\item The map $\Phi_n^\mathcal{F}$ induces the identity natural transformation  on $\HR^\mathcal{F}_{n}$.
\item For any $n\geq 1$, there exists non-thin elements $x$ in $\mathcal{F}_{n-1}(\C)$
as soon as $\C_n\neq \emptyset$. Indeed, if $u\in \C_n$,  $\ev\ \square_n^\mathcal{F}(u)=u$,
therefore $\square_n^\mathcal{F}(u)$ is a non-thin element of $\mathcal{F}_{n-1}(\C)$.
\end{itemize}

We end this section by some general facts about regular cuts.

\bp\label{Phi_functoriel} Let $f$ be a morphism of cuts from $\mathcal{F}$ to $\mathcal{G}$. Suppose that $\mathcal{F}$ and
$\mathcal{G}$ are regular. Then $\Phi^\mathcal{G} \circ f=f\circ \Phi^\mathcal{F}$ as natural
transformation from $\mathcal{F}$ to $\mathcal{G}$. In other terms, the following diagram
is commutative :
\[ \xymatrix{\mathcal{F} \fr{f} \fd{\Phi^\mathcal{F}}& \mathcal{G} \ar@{->}[d]^{\Phi^\mathcal{G}} \\
\mathcal{F} \fr{f} & \mathcal{G}} \]
\ep

\bpf Let $n\geq 0$ and let $P(n)$ be the property : ``for any $\omega$-category
$\C$ and any $x\in tr^{n}\P\C$, then $f\square_{n+1}^\mathcal{F}(x)= \square_{n+1}^\mathcal{G} x$.''

One has $\Phi_1^\mathcal{F}:=\id_{\mathcal{F}_0}$, $\Phi_1^\mathcal{G}:=\id_{G_0}$ and necessarily
$f_0=I\!d$ by definition of a morphism of cuts.
Therefore $P(0)$ holds. Now suppose
$P(n)$ proved for some $n\geq 0$.  One has
$\ev f \square_{n+2}^\mathcal{F} = \ev \square_{n+2}^\mathcal{F} = Id_{tr^{n+1}\P}$ since
$f$ is a morphism of cuts and
\begin{alignat*}{2}
f \square_{n+2}^\mathcal{F} i_{n+1}  &= f (\epsilon_n.\square_{n+1}^\mathcal{F}) i_{n+1} &&\\
&= f \epsilon_n \square_{n+1}^\mathcal{F}  &&\hbox{ by definition of $\epsilon_n.\square_{n+1}^\mathcal{F}$}\\
&=\epsilon_n f \square_{n+1}^\mathcal{F} &&\hbox{ since $f$ morphism of simplicial sets}\\
&= \epsilon_n \square_{n+1}^\mathcal{G} && \hbox{ by induction hypothesis}
\end{alignat*}
Therefore the natural transformation $f \square_{n+2}^\mathcal{F}$ from $tr^{n+1}\P$
to $\mathcal{G}_{n+1}$ can be identified with $\epsilon_n.\square_{n+1}^\mathcal{G}$ which is
precisely $\square_{n+2}^\mathcal{G}$. Therefore $P(n+1)$ is proved.

At last, if $x\in \mathcal{F}_n(\C)$, then
\begin{alignat*}{2}
\Phi^\mathcal{G} f(x)&=\square_{n+1}^\mathcal{G}\ev f(x)&&\ \hbox{ by definition of folding operators}\\
&=\square_{n+1}^\mathcal{G} \ev (x) &&\ \hbox{ since $f$ preserves the evaluation map}\\
&= f\square_{n+1}^\mathcal{F} \ev (x) &&\ \hbox{ since $P(n)$ holds}\\
&= f \Phi^\mathcal{F}(x)&&\ \hbox{  by definition of folding operators}
\end{alignat*}

\epf

\bp\label{canonical_position} If $u$ is a $(n+1)$-morphism of $\C$ with $n\geq 1$,
then $\square_{n+1}^\mathcal{F} u$ is an homotopy within the
simplicial set $\mathcal{F}(\C)$ between $\square_n^\mathcal{F} s_nu$
and $\square_n^\mathcal{F} t_nu$. \ep

\bpf The natural map $\ev\ \de_i\square_{n+1}^\mathcal{F}$ for $0\leq i\leq n$
from $tr^{n}\P$ to $tr^{n-1}\P$ is of the form
$d_{m_i}^{\alpha_i}$ for $m_i\leq n$ with $m_i\leq n-1$ for $0\leq i\leq n-2$
and  $\left\{\ev \de_{n-1}\square_{n+1}^\mathcal{F}, \ev \de_{n}\square_{n+1}^\mathcal{F}\right\}=\left\{s_{n},t_{n}\right\}$. Therefore for $0\leq i\leq n-2$,
$\de_{i}\square_{n+1}^\mathcal{F}=\de_{i}\square_{n+1}^\mathcal{F}s_n= \de_{i}\square_{n+1}^\mathcal{F}t_n$ by rule~\ref{regular3.02} of Definition~\ref{def_regular}. And
by construction of $\square_{n+1}^\mathcal{F}$, one obtains
$\de_{i}\square_{n+1}^\mathcal{F}=\epsilon_{n-2}\de_{i}\square_{n}^\mathcal{F}s_n=\epsilon_{n-2}\de_{i}\square_{n}^\mathcal{F}t_n$.
\epf

\begin{cor}\label{diff_reduite}
If  $x\in \CR^\mathcal{F}_{n+1}(\C)$, then
$\de x=\de\square_{n+1}^\mathcal{F}x=\square_{n}^\mathcal{F} s_nx-\square_{n}^\mathcal{F}t_nx$ in
$\CR^\mathcal{F}_n(\C)$.  In other terms, the differential map from
$\CR^\mathcal{F}_{n+1}(\C)$ to $\CR^\mathcal{F}_{n}(\C)$ with $n\geq 1$ is induced by the
map $s_n-t_n$.
\end{cor}

\section{The cuts of branching and merging ner\-ves}\label{example_corner}

We see now that the corner nerves $\mathcal{N}^\eta$ defined in
\cite{Gau} are two examples of regular cuts
with the correspondence
$\square_n^\eta:=\square_n^{\mathcal{N}^\eta}$,
$\Phi_n^\eta:=\Phi_n^{\mathcal{N}^\eta}$,
$H_n^\eta:=H_n^{\mathcal{N}^\eta}$, $\HR_n^\eta:=\HR_n^{\mathcal{N}^\eta}$
and $\ev(x)=x(0_{dim(x)})$.

Let us first recall the construction of the free $\omega$-category
$I^n$ generated by the faces of the $n$-cube. The faces of the
$n$-cube are labeled by the words of length $n$ in the alphabet
$\{-,0,+\}$, one word corresponding to the barycenter of one
face. We take the convention that $00\dots 0 \hbox{ ($n$
times)}=:0_n$ corresponds to its interior and that $-_n$ (resp.
$+_n$) corresponds to its initial state $--\dots -\hbox{ ($n$
times)}$ (resp. to its final state $++\dots +\hbox{ ($n$
times)}$).  If $x$ is a face of the $n$-cube, 
let $R(x)$ be the set of faces of
$x$. If $X$ is a set of faces, then let $R(X)=\bigcup_{x\in
X}R(x)$. Notice that $R(X\cup Y)=R(X)\cup R(Y)$ and that
$R(\{x\})=R(x)$. Then $I^n$ is the free $\omega$-category
generated by the $R(x)$ with the rules

\begin{enumerate}
\item  For $x$ $p$-dimensional with $p\geq 1$,
$$s_{p-1}(R(x))=R(s_x)$$
and
$$t_{p-1}(R(x))=R(t_x)$$
where $s_x$ and
$t_x$ are the sets of faces defined below.
\item If $X$ and $Y$ are two elements of $I^n$ such that
$t_p(X)=s_p(Y)$ for some $p$, then $X\cup Y$ belongs to $I^n$
 and $X\cup Y=X *_p Y$.
\end{enumerate}

The set $s_x$ is the set of subfaces of the faces obtained by
replacing the $i$-th zero of $x$ by $(-)^i$, and the set $t_x$ is the
set of subfaces of the faces obtained by replacing the $i$-th zero of
$x$ by $(-)^{i+1}$. For example,
$s_{0+00}=\{\hbox{-+00},\hbox{0++0},\hbox{0+0-}\}$ and
$t_{0+00}=\{\hbox{++00},\hbox{0+-0},\hbox{0+0+}\}$. Figure~\ref{I3}
represents the free $\omega$-category generated by the $3$-cube.

The branching and merging nerves are dual from each other.
We set
\beas
&&\omega Cat(I^{n+1},\C)^\eta:=\{x\in
\omega Cat(I^{n+1},\C),d_0^\eta(u)=\eta_{n+1} \\
&&\hbox{ and }dim(u)=1 \Longrightarrow dim(x(u))=1\}
\eeas
where $\eta\in\{-,+\}$ and where $\eta_{n+1}$ is the initial state
(resp. final state) of $I^{n+1}$ if $\eta=-$ (resp. $\eta=+$). For all
$(i,n)$ such that  $0\leq i\leq n$, the face maps $\de_i$ from
$\omega Cat(I^{n+1},\C)^\eta$ to $\omega Cat(I^{n},\C)^\eta$
are the arrows $\de^\eta_{i+1}$ defined by
\[\de^\eta_{i+1}(x)(k_1\dots k_{n+1})=x(k_1\dots [\eta]_{i+1}\dots k_{n+1})\]
and the degeneracy maps $\epsilon_i$ from $\omega Cat(I^{n},\C)^\eta$
to $\omega Cat(I^{n+1},\C)^\eta$
are the arrows $\Gamma^\eta_{i+1}$ defined by setting
\begin{eqnarray*}
&& \Gamma_i^-(x)(k_1\dots k_n):=x(k_1\dots \max(k_i,k_{i+1})\dots k_n)\\
&&\Gamma_i^+(x)(k_1\dots k_n):=x(k_1\dots \min(k_i,k_{i+1})\dots k_n)
\end{eqnarray*}
with the order $-<0<+$.

\begin{propdef}\cite{Gau}\label{def-coin} Let $\C$ be an
  $\omega$-category. The $\N$-graded set $\mathcal{N}^\eta(\C)$ together
with the convention $\mathcal{N}^\eta_{-1}(\C)=\C_0$, endowed with
the maps $\de_i$ and $\epsilon_i$ above defined with moreover
$\de_{-1}=s_0$ (resp. $\de_{-1}=t_0$) if $\eta=-$ (resp. $\eta=+$) and
with $\ev(x)=x(0_n)$ for $x\in \omega Cat(I^{n},\C)$ is a simplicial
cut. It is called the  \textit{$\eta$-corner
simplicial nerve} $\mathcal{N}^\eta$ of $\C$.
\end{propdef}

Set
$H_{n+1}^{\eta}(\C):=H_n(\mathcal{N}^{\eta}(\C))$
for $n\geq -1$.  These
homology theories are called \textit{branching} and \textit{merging homology}
respectively and are exactly the same homology theories as that
defined in \cite{Gau} and studied in \cite{Coin}.

And we have

\bth\cite{Coin} The simplicial cut $\mathcal{N}^\eta$ is regular. The associated
folding operator $\square_n^{\mathcal{N}^\eta}$ coincides with the
operator $\square_n^\eta$ defined in \cite{Coin}. And therefore the
associated homology theory $\HR_n^{\mathcal{N}^\eta}$ coincide with
the reduced corner homology $\HR_n^\eta$ defined in \cite{Coin}. \eth

It is useful for the sequel to remind some important properties of the
folding operators associated to corner nerves.

\bth\label{caracterisation}\cite{Coin} Let
$\C$ be an $\omega$-category. Let $x$ be an element of
$\mathcal{N}^-_n(\C)$. Then the following two conditions are equivalent :
\begin{enumerate}
\item the equality $x=\Phi_n^-(x)$ holds
\item for $1\leq i \leq n$, one
has $\ev \de_i^+ x= \de_i^+ x(0_n)$ is $0$-dimensional
and for $1\leq i \leq n-2$, one
has $\de_i^- x\in Im(\Gamma_{n-2}^-\dots \Gamma_i^-)$.
\end{enumerate}
\eth

Another operator coming from \cite{Coin} which matters for this paper
is the operator $\theta_i^-$.

\bd Let $x\in \mathcal{N}_n^-(\C)$ for some $\C$ such that for any
$1\leq j\leq n+1$, $\de_j^+x$ is $0$-dimensional. Then $x$ is called a
\textit{negative} element of the branching nerve.
\ed

\bth\label{decomposition_theta}
Let $n\geq 2$. There exists natural transformations
$$\theta_1^-,\dots,\theta_{n-1}^-$$ from $\mathcal{N}_n^-$ to itself
satisfying the following properties :
\begin{enumerate}
\item\label{th1} If $x$ is a negative element of $\mathcal{N}_n^-(\C)$, then for
  any $1\leq i\leq n-1$, $\theta_i^-x$ is a negative element as well.
\item\label{th2} If $x$ is a negative element of $\mathcal{N}_n^-(\C)$, then for
  any $1\leq i\leq n-1$, there exists a thin negative element $y_i$ of
  $\mathcal{N}_{n+1}^-(\C)$ such that $\de^-y_i-x$ is a linear
  combination of thin negative elements.
\item\label{th3} There exists a composite of $\theta_1^-,\dots,\theta_{n-1}^-$ which
coincides with the negative folding operators on negative elements of
$\mathcal{N}_n^-$.
\end{enumerate}
\eth

\bpf[Sketch of proof] Consider  the $\theta_1^-,\dots,\theta_{n-1}^-$ of \cite{Coin}.
One has
\beas
&&\de_j^+ \theta_i^- =\left\{
\begin{array}{c} \theta_{i-1}^- \de_j^+ \hbox{ if $j<i$}\\
\theta_{i}^- \de_j^+
\hbox{ if $j>i+2$}\end{array}\right.\\
&& \de_i^+ \theta_i^- = \ {}^v\psi_i^- \de_i^+ \\
&& \de_{i+1}^+ \theta_i^- = \epsilon_{i+1} \de_{i+1}^+ \de_i^- +_i \epsilon_{i+1} \de_{i+1}^+ \de_{i+1}^+\\
&& \de_{i+2}^+ \theta_i^- = \ {}^v\psi_i^+ \de_{i+2}^+
\eeas
where, for the last formula, $\ {}^v\psi_i^\pm $ are other operators which is
not important to explicitly define  here : the only important thing is that $\de_i^+ \theta_i^-$
remains $0$-dimensional if the argument is $0$-dimensional. Hence property~\ref{th1}.
As for property~\ref{th2}, it is enough to check  it for $i=1$. And in this
case, $y$ is a  thin $4$-cube satisfying
\beas
&&\de_1^+ y= \ {}^v\psi_2^- \Gamma_1^- \de_1^+ x\\
&&\de_2^+ y= \Gamma_2^- \de_2^+ x\\
&& \de_3^+ y= \epsilon_3(\Gamma_1^- \de_2^+ \de_1^- x +_1 \epsilon_2\de_2^+ \de_2^+ x)\\
&&\de_4^+ y=\ {}^v\psi_2^+ \Gamma_2^- \de_3^+ x
\eeas
Once again, we refer to \cite{Coin} for the precise definition of
the  operators involved in the above formulas.
The only thing that matters here is the dimension of
$\de_i^+y$.

By \cite{Coin}, we know that $\Phi^-=\Theta \circ \Psi$ when $\Theta$ is a
composite of $\theta_i^-$ and such that for $x$ negative, $\Psi x=x$.
Hence property~\ref{th3}.
\epf

The graded set $(\omega Cat(I^n,\C))_{n\geq 0}$ endowed with the operations
$\de_i^\pm$ above defined and by the maps
$\epsilon_i(x)(k_1\dots k_{n+1})=x(k_1\dots \widehat{k_i}\dots k_{n+1})$ for
$x\in \omega Cat(I^n,\C)$ and $1\leq i\leq n+1$
is a cubical set and is usually known as the \textit{cubical singular
  nerve} of $\C$ \cite{Brown_cube}. The use of the same notation
$\epsilon_i$ for the degeneracy maps of the cubical singular nerve and
the degeneracy maps of the three simplicial nerves appearing in this
paper is very confusing. Fortunately, we will not need the degeneracy
maps of the cubical singular nerve in this work except for Theorem~\ref{decomposition_theta} right above.

\section{The globular cut}\label{glob_cut}

The most direct way of constructing a cut of $\omega$-categories consists
of using the composite of both functors $\P:\C\mapsto \P\C$ and
$\mathcal{N}$ where $\mathcal{N}$ is the simplicial nerve functor
defined by Street~\footnote{Of course, the functor $\mathcal{N}$ can
  be viewed as a functor from $\omega Cat_1$ to
  $Sets^{\Delta^{op}}_+$, but a ``good'' cut should not be extendable
  to a functor from $\omega Cat$ to $Sets^{\Delta^{op}}_+$.}.

Let us start this section by recalling the construction of the
free $\omega$-catego\-ry $\Delta^n$ generated by the faces of the
$n$-simplex. The faces of the $n$-simplex are labeled by the
strictly increasing sequences of elements of $\{0,1,\dots,n\}$.
The length of a sequence is equal to the dimension of the
corresponding face plus one. If $x$ is a face of the $n$-simplex, its
subfaces are all increasing sequences of $\{0,1,\dots,n\}$
included in $x$. If $x$ is a face of the $n$-simplex, 
let $R(x)$ be the set of faces
of $x$. If $X$ is a set of faces, then let $R(X)=\bigcup_{x\in
X}R(x)$. Notice that $R(X\cup Y)=R(X)\cup R(Y)$ and that
$R(\{x\})=R(x)$. Then $\Delta^n$ is the free $\omega$-category
generated by the $R(x)$ with the rules

\begin{enumerate}
\item  For $x$ $p$-dimensional with $p\geq 1$,
$$s_{p-1}(R(x))=R(s_x)$$ and $$t_{p-1}(R(x))=R(t_x)$$
where $s_x$ and
$t_x$ are the sets of faces defined below.
\item If $X$ and $Y$ are two elements of $\Delta^n$ such that
$t_p(X)=s_p(Y)$ for some $p$, then $X\cup Y$ belongs to $\Delta^n$
and $X\cup Y=X *_p Y$.
\end{enumerate}

where $s_x$ (resp. $t_x$) is the set of subfaces of $x$ obtained by
removing one element in odd position (resp. in even position). For
instance, $s_{(04589)}=\{(4589),(0489),(0458)\}$ and
$t_{(04589)}=\{(0589),(0459)\}$.

Sometimes we will write (for instance) $(0<4<5<8<9)$ instead of simply
$(04589)$. Figure~\ref{2simpl} gives the example of the $2$-simplex.

Let $x\in \omega Cat(\Delta^n,\C)$.
Then consider the labeling of
the faces of respectively $\Delta^{n+1}$ and $\Delta^{n-1}$ defined by :
\begin{itemize}
\item $\epsilon_i(x)(\sigma_0<\dots<\sigma_{r})=x(\sigma_0<\dots<\sigma_{k-1}<\sigma_{k}-1<\dots<\sigma_{r}-1)$
 if $\sigma_{k-1}<i$ and $\sigma_{k}>i$.
\item $x(\sigma_0<\dots<\sigma_{k-1}<i<\sigma_{k+1}-1<\dots<\sigma_{r}-1)$
if $\sigma_{k-1}<i$, $\sigma_{k}=i$ and $\sigma_{k+1}>i+1$.
\item $x(\sigma_0<\dots<\sigma_{k-1}<i<\sigma_{k+2}-1<\dots<\sigma_{r}-1)$
if $\sigma_{k-1}<i$, $\sigma_{k}=i$ and $\sigma_{k+1}=i+1$.
\end{itemize}
and
\[
\de_i(x)(\sigma_0<\dots<\sigma_{s})=x(\sigma_0<\dots<\sigma_{k-1}<\sigma_{k}+1<\dots<\sigma_{s}+1)
\]
where $\sigma_k,\dots,\sigma_{s}\geq i$ and $\sigma_{k-1}<i$.

It can be checked that $\epsilon_i(x)$ (resp. $\de_i(x)$) are
$\omega$-functors from $\Delta^{n+1}$ (resp. $\Delta^{n-1}$) to $\C$
\cite{oriental}. By construction, the map $[n]\mapsto \Delta^n$ induces
then a functor from the well-known category $\Delta$ whose associated
presheaves are the simplicial sets to $\omega Cat$. Therefore
$\mathcal{N}(\C)=(\omega Cat(\Delta^*,\C),\de_i,\epsilon_i)$ is a
simplicial set which is called the \textit{simplicial nerve} of $\C$.

\bd\label{def_glob_nerve} The globular cut $\mathcal{N}^{gl}$ (or the globular nerve)
is the functor from $\omega Cat_1$ to $Sets^{\Delta^{op}}_+$ defined
by $\mathcal{N}^{gl}_n(\C)=\omega Cat(\Delta^n,\P\C)$ for $n\geq 0$
and with $\mathcal{N}^{gl}_{-1}(\C)=\C_0\p\C_0$, and endowed with the
augmentation map $\de_{-1}$ from $\mathcal{N}^{gl}_0(\C)=\C_1$ to
$\mathcal{N}^{gl}_{-1}(\C)=\C_0\p \C_0$ defined by
$\de_{-1}x=(s_0x,t_0 x)$. The evaluation map $\ev$ is defined by
$\ev(x)=x((0\dots n))$ for $x\in \omega Cat(\Delta^n,\P\C)$.  The
homology theory $H_n^{gl}:=H_n^{\mathcal{N}^{gl}}$ is called the
globular homology and $\HR_n^{gl}:=\HR_n^{\mathcal{N}^{gl}}$ the
reduced globular homology.
\ed

\begin{figure}
\begin{center}
\includegraphics[width=8cm]{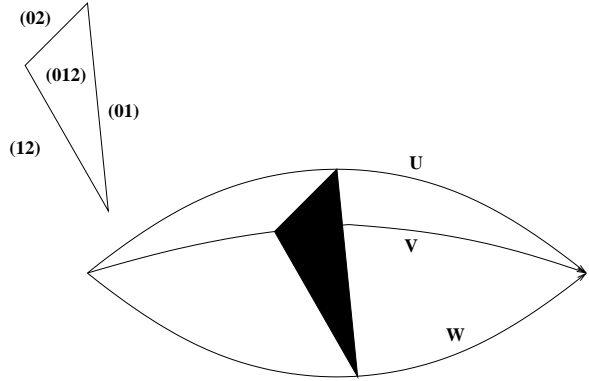}
\end{center}
\caption{Globular $2$-simplex}
\label{2simplglob}
\end{figure}

Geometrically, the elements of $\mathcal{N}^{gl}_n(\C)$ are full
$(n+1)$-globes. Figure~\ref{2simplglob} depicts a $2$-simplex in
the globular nerve. The simplexes seen by the globular cut are
intuitively transverse to the execution paths, as well as those
of corner nerves. Hence the terminology of cuts.

Here is now the new definition of the globular homology of a 
globular 
$\omega$-category $\C$ :

\bd\label{def_correcte}
Let $\C$ be a non-contracting $\omega$-category. We set
\[H_{n+1}^{gl}(\C):=H_n(\mathcal{N}^{gl}(\C))\] for $n\geq
-1$ and this homology theory is called the globular homology of $\C$.
\ed

\section{Associating to any globe its corners}\label{morphisms}

\begin{figure}
\begin{center}
\subfigure[A $2$-globular simplex $\widetilde{X}$]{\label{2sim}
\xymatrix{&& {w}\\&&\\{u}\ar@{->}[rruu]^{C}\ar@{->}[rr]_{A}&&{v}\ar@{->}[uu]_{B}
\ar@2{->}[lu]|{X}}
}
\end{center}
\begin{center}
\subfigure[The $2$-simplex $h^-(\widetilde{X})$]{\label{h2sim}
\xymatrix{
&\ar@{->}[rr]^{t_0u}& &\ar@{->}[rd]^{t_0u}&&&
&\ar@{->}[dr]|{t_0u}\ar@{->}[rr]^{t_0u}&&\ar@{->}[dr]^{t_0u}&\\
\ar@{->}[ru]^{w}\ar@{->}[rr]|{v}\ar@{->}[rd]_{u}&&\ar@{->}[ru]|{t_0u}\ar@{->}[dr]|{t_0u}\ff{lu}{B}&&
\ar@3{->}[rr]^{X}&&
\ar@{->}[ur]^{w}\ar@{->}[dr]_{u}&&\ff{ur}{t_0u}\ar@{->}[rr]|{t_0u}&&\\
&\ar@{->}[rr]_{t_0u}\ff{ru}{A}&&\ar@{->}[ru]_{t_0u}\ff{uu}{t_0u}&&&
&\ar@{->}[ru]|{t_0u}\ar@{->}[rr]_{t_0u}\ff{uu}{C}&&\ff{ul}{t_0u}\ar@{->}[ru]_{t_0u}\\
}
}
\end{center}
\caption{Illustration of $h^-$}
\label{bleu}
\end{figure}
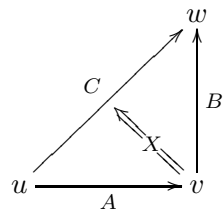
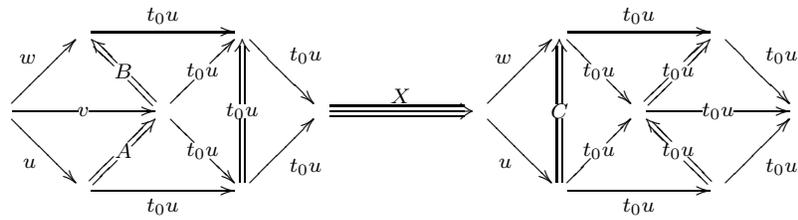

The purpose of the rest of the paper is to justify that
Definition~\ref{def_correcte} is the right definition. This is
not a mathematical statement of course ! We follow the order of
the remarks at the very end of Section~\ref{intro} which explain
what kind of conditions the globular homology must fulfill. So we
have first to construct $h^-$ and $h^+$ and we must verify that
geometrically, in homology, $h^-$ and $h^+$ do what we expect to
find.  In fact, we refer to \cite{ConcuToAlgTopo} for intuitive
explanations of $h^-$ and $h^+$. We only recall here
Figure~\ref{bleu} as an illustration and care only about the
construction of $h^-$.

\bth\label{hmoins} Let $\alpha\in\{-,+\}$.  There exists one and only
one morphism of cuts $h^\alpha$ from $\mathcal{N}^{gl}$ to
$\mathcal{N}^{\alpha}$. Moreover, for any non-contracting 
$\omega$-category $\C$, both morphisms $h^\alpha$ from
$\mathcal{N}^{gl}(\C)$ to $\mathcal{N}^{\alpha}(\C)$ are injective.
\eth

The rest of the section is devoted to the proof of
Theorem~\ref{hmoins}. The following sequence of propositions establishes the
existence of $h^-$. The term $\Cube^n$ denotes the set of faces of
the $n$-cube, as described in Section~\ref{example_corner}.

We briefly recall how filling shells in the cubical singular nerve. This
technical tool already appears in \cite{Brown_cube} for
$\omega$-groupoids and in \cite{phd-Al-Agl} for $\omega$-cate\-go\-ries. A
particular case can be found in \cite{Gau}.

\bd A \textit{$n$-shell} in the cubical singular nerve is a family of
$2(n+1)$ elements $x_i^\pm $ of $\omega Cat(I^n,\C)^-$ such that
$\de_i^\alpha x_j^\beta= \de_{j-1}^\beta x_i^\alpha$ for $1\leq
i< j\leq n+1$ and $\alpha,\beta\in\{-,+\}$.  \ed

If $x_i^\pm $ is a $n$-shell, then it induces a labeling $x$ on the
set of faces of dimension 
at most $n$ of the $(n+1)$-cube in the following manner :
let $k_1\dots k_{n+1}$ be a face of dimension at most $n$ ; then there
exists $i$ such that $k_i\neq 0$ ; then let $x(k_1\dots k_{n+1}):=
x_i(k_1\dots \widehat{k_i} \dots k_{n+1})$. The axiom satisfied by an
$n$-shell ensures the coherence of the definition.

\begin{propdef}\label{filling}
Let $x_i^\pm $ be an $(n-1)$-shell with $n\geq 1$.
\begin{itemize}
\item The labeling of the faces of dimension at most $(n-1)$ of $I^n$ 
defined  by  $x_i^\pm $ always induces
an $\omega$-functor and only one from $I^n\backslash \{R(0_n)\}$
to $\C$. Denote it by $x$.
\item The $n$-shell $(x_i^\pm )$ is said \textit{fillable}
  if there exists a morphism $u$ of $\C$ such that
$s_{n-1}u= x\left( s_{n-1} R(0_n) \right)$ and $t_{n-1}u= x\left( t_{n-1} R(0_n) \right)$. In this case, there
exists a unique $\omega$-functor $x$ from $I^n$ to $\C$ such that
$\de_i^\pm x= x_i^\pm$ for $1\leq i\leq n$ and $x(0_n)=u$.
\end{itemize}
\end{propdef}

\bpf Using the
freeness of $I^n$, the construction in the proof of \cite{Gau}
Proposition 5.1 yields the $\omega$-functor $x$ from $I^n\backslash \{R(0_n)\}$
to $\C$. The hypotheses stated in \cite{Gau}  were too strong indeed.  If moreover the shell
is fillable in the above sense, one concludes still as in the proof of \cite{Gau}
Proposition 5.1.
\epf

Now we can construct $h^-$.

\bth\label{simplicial_h} Let $x$ be an $n$-simplex of the globular simplicial nerve of
$\C$. Then the map $h_n^-(x)$ from $\Cube^{n+1}$ to $\C$ defined by
\begin{enumerate}
\item\label{re1}  $+\in \{k_1\dots k_{n+1}\}$ implies
$h_n^-(x)(k_1\dots k_{n+1})=t_0 x((0))$ (notice that $(0)$ is the final state of
$\Delta^n$)
\item\label{re2} $\{k_1,\dots ,k_{n+1}\}\subset \{-,0\}$ and $$\{k_1,\dots, k_{n+1}\}\cap \{0\}=
\{k_{\sigma_0+1},\dots,k_{\sigma_r+1}\}$$ with $\sigma_0<\dots<\sigma_r$ implies
$h_n^-(x)(k_1\dots k_{n+1})=x((\sigma_0\dots\sigma_r))$
\item\label{re3} $h_n^-(x)(-_{n+1})=s_0 x((n))$ (notice that $(n)$ is the initial state of
$\Delta^n$)
\end{enumerate}
yields an $\omega$-functor from $I^{n+1}$ to $\C$. Moreover, $h^-$
induces a morphism of simplicial sets from the globular nerve of $\C$
to its negative corner nerve. And the map from
$\mathcal{N}^{gl}_{-1}(\C)$ to $\mathcal{N}^{-}_{-1}(\C)$ defined by
$(x,y)\mapsto x$ extends the previous morphism to the corresponding
augmented simplicial nerves. Moreover for $n\geq 0$, $h_n^-$ is a
one-to-one map and the image of $h_n^-$ contains exactly all cubes $x$
of the negative corner nerve such that as soon as $\de_i^+x$ exists,
then it is $0$-dimensional. \eth

There
is no ambiguity to set $h^-(x)=h_n^-(x)$ if $x$ is an $n$-simplex of
the globular cut.

In the sequel, in order to make easier the reading of  the calculations, we suppose
that an expression like $(\sigma_0<\sigma_j\leq \widehat{k} <
\sigma_{j+1}<...<\sigma_r)$ is the same thing as $(\sigma_0<\sigma_j<
\sigma_{j+1}<...<\sigma_r)$ in $\Delta^*$ but with an additional
information given within the calculation itself : here that
$\sigma_j\leq \widehat{k}< \sigma_{j+1}$ holds.

\begin{proof}
One proves by induction on $n$ the following property $P(n)$ : `` For
any $n$-simplex $x$ of the globular simplicial nerve of any
$\omega$-category $\C$, the map $h^-(x)$ from $\Cube^{n+1}$ to $\C$
induces an $\omega$-functor and moreover an element of $\omega
Cat(I^{n+1},\C)^-$.''

Let $x$ be a $0$-simplex of the globular nerve of $\C$. Then $x$ is an
$\omega$-functor from $\Delta^0$ to $\P\C$, and therefore it can be
identified with the $1$-morphism $x((0))$ of $\C$. Therefore
\begin{alignat*}{2}
&h^-(x)(0)=x((0)) &&  \ \hbox{by rule~\ref{re2}}\\
& h^-(x)(+)=t_0 x((0))&&  \ \hbox{by rule~\ref{re1}}\\
& h^-(x)(-)=s_0 x((0)) && \ \hbox{by rule~\ref{re3}}
\end{alignat*}
Therefore $P(0)$ is proved.

Now suppose that $P(n)$ is proved for $n\geq
0$. Let $x$ be a $(n+1)$-simplex of the globular simplicial nerve of
some $\omega$-category $\C$.  If $+\in\{k_1,\dots, k_{n+1}\}$, then
\begin{alignat*}{2}
&\de_i^-(h^-(x))(k_1\dots k_{n+1})&&\ \\
&=h^-(x)(k_1\dots k_{i-1} - k_i\dots k_{n+1})
&&\ \hbox{ by definition of $\de_i^-$ for $1\leq i\leq n+2$}\\
&=t_0 x((0))&&\ \hbox{ by rule~\ref{re1}}\\
&=h^-(\de_{i-1}x)(k_1\dots k_{n+1})&&\ \hbox{ again by rule~\ref{re1}}
\end{alignat*}
If $+\notin\{k_1,\dots, k_{n+1}\}$, i.e. if
$\{k_1,\dots, k_{n+1}\}\subset \{-,0\}$, set
\[\{k_1,\dots, k_{n+1}\}\cap \{0\}= \{k_{\sigma_0+1},\dots,k_{\sigma_r+1}\}\]
with $\sigma_0<\dots<\sigma_r$. For a given
$i$ such that $1\leq i\leq n+2$, set $$w_1\dots w_{n+2}=k_1\dots k_{i-1} - k_i\dots k_{n+1}$$
as word. Then  let $$\{w_1,\dots, w_{n+2}\}\cap \{0\}=\{w_{\tau_0+1},\dots,w_{\tau_r+1}\}$$
with $\tau_0<\dots<\tau_r$. The relation between the sequence of $\sigma_j$ and the
sequence of $\tau_j$ is as follows :
\beas
&&\sigma_j+1\leq i-1\Longrightarrow \sigma_j=\tau_j\\
&&\sigma_j+1\geq i \Longrightarrow \sigma_j+1=\tau_j
\eeas
And we have
\begin{alignat*}{2}
&\de_i^-(h^-(x))(k_1\dots k_{n+1})&&\ \\
&=h^-(x)(k_1\dots k_{i-1} - k_i\dots k_{n+1})\ \hbox{ by definition of $\de_i^-$}
&&\\
&= x((\tau_0\dots\tau_r))\ \hbox{ by rule~\ref{re2}}&&\\
&= x((\sigma_0<\dots <\sigma_{j_0}\leq \widehat{i-2}<\widehat{i-1}< \sigma_{j_0+1}+1<\dots<\sigma_r+1))&&\\
&=(\de_{i-1}x)((\sigma_0\dots\sigma_r))\  \hbox{ by definition of $\de_{i-1}$}&&\\
&=h^-(\de_{i-1}x)(k_1\dots k_{n+1})\ \hbox{ by rule~\ref{re2}} &&
\end{alignat*}
Therefore $\de_i^-(h^-(x))=h^-(\de_{i-1}x)$.  And by
rule~\ref{re1}, $\de_i^+ (h^-(x))$ is the constant
$\omega$-functor from $\Cube^{n+1}$ to $\C$ which sends any face of
$I^{n+1}$ on $t_0 x((0))$. Therefore $(\de_i^\pm (h^-(x)))_{1\leq i\leq n+1}$
is a $(n+1)$-shell in the cubical
nerve of $\C$ which is fillable. By Proposition~\ref{filling}, the
labeling $h^-(x)$ of $\Cube^{n+2}$ induces an
$\omega$-functor from $I^{n+2}$ to $\C$ and $P(n+1)$ is proved.

By construction, the equality $\de_i^- (h^-(x))=h^-(\de_{i-1}x)$ holds for any
$n$-simplex $x$ of the globular nerve and for $1\leq i\leq n+1$. It remains
to check that for such a simplex $x$, $\Gamma_i^-(h^-(x))= h^-(\epsilon_{i-1}x)$
for $i\leq 1\leq n+1$. Consider a face $k_1\dots k_{n+2}$ of the $(n+2)$-cube.
If $+\in\{k_1,\dots, k_{n+2}\}$, then
\begin{alignat*}{2}
& \Gamma_i^-(h^-(x))(k_1\dots k_{n+2}) &&\ \\
&= h^-(x)(k_1\dots max(k_i,k_{i+1})\dots k_{n+2})&&\ \hbox{ by definition of $\Gamma_i^-$}\\
&= t_0 x((0))&&\ \hbox{ by rule~\ref{re1}}\\
&= h^-(\epsilon_{i-1}x)(k_1\dots k_{n+2})&&\ \hbox{ again by rule~\ref{re1}}
\end{alignat*}
If $+\notin\{k_1,\dots, k_{n+2}\}$, i.e. if
$\{k_1,\dots, k_{n+2}\}\subset \{-,0\}$, set
\[\{k_1,\dots, k_{n+2}\}\cap \{0\}= \{k_{\sigma_0+1},\dots,k_{\sigma_r+1}\}\]
with $\sigma_0<\dots<\sigma_r$. For a given
$i$ such that $1\leq i\leq n+1$,
\[\{k_1,\dots, max(k_i,k_{i+1}),\dots, k_{n+2}\}\subset \{-,0\}\] and set
$w_1\dots w_{n+1}=k_1\dots max(k_i,k_{i+1})\dots k_{n+2}$
as word. Then let
\[\{w_1,\dots, w_{n+1}\}\cap \{0\}= \{w_{\tau_0+1},\dots,w_{\tau_s+1}\}\]
with $\tau_0<\dots<\tau_s$. One has to calculate
\begin{alignat*}{2}
& \Gamma_i^-(h^-(x))(k_1\dots k_{n+2}) &&\ \\
&= h^-(x)(k_1\dots max(k_i,k_{i+1})\dots k_{n+2})&&\ \hbox{ by definition of $\Gamma_i^-$}\\
&= x((\tau_0\dots \tau_s)) &&\ \hbox{ by definition of $h^-$}
\end{alignat*}
for some $1\leq i \leq n+2$.

The situation can be decomposed in three mutually exclusive cases\thinspace:

\begin{enumerate}
\item $k_i=k_{i+1}=0$. In this case, there exists a unique $j_0$ such that
$\sigma_{j_0}+1=i$, $s=r-1$ and
\beas
&&\sigma_j+1\leq i-1 \Longrightarrow \sigma_j=\tau_j \hbox{ (in this case, $j<j_0$)}\\
&&\tau_{j_0}+1=i=\sigma_{j_0}+1\\
&&\sigma_j+1\geq i+2 \Longrightarrow \sigma_j- 1=\tau_{j-1}\hbox{ (in this case, $j>j_0+1$)}
\eeas
Then $\sigma_{j_0+2}\geq i+1$ and
\begin{alignat*}{2}
& x((\tau_0\dots \tau_s))&&\ \\
&= x((\sigma_0<\dots <\sigma_{j_0}=\widehat{i-1}<\sigma_{j_0+2}-1<\dots<\sigma_{s+1} -1))&&\ \\
&= (\epsilon_{i-1} x)(\sigma_0\dots\sigma_{j_0}\sigma_{j_0+1}\sigma_{j_0+2}\dots\sigma_{s+1})\ \hbox{ by definition of $\epsilon_i$} &&\\
& \ \hbox{  and since $\sigma_{j_0+1}=i$} &&\\
&= (h^-(\epsilon_{i-1} x))(k_1\dots k_{n+2})\ \hbox{ by definition of $h^-$}&&
\end{alignat*}
\item $k_i= k_{i+1}=-$. In this case, $s=r$ and
\beas
&&\sigma_j+1\leq i-1 \Longrightarrow \sigma_j=\tau_j\\
&&\sigma_j+1\geq i+2 \Longrightarrow \sigma_j- 1=\tau_{j}
\eeas
Then for some $k$,
\begin{alignat*}{2}
& x((\tau_0\dots \tau_s))&&\ \\
&= x((\sigma_0<\dots <\sigma_{k}<\widehat{i-1}<\sigma_{k+1}-1<\dots<\sigma_{r}-1))&&\\
&= (\epsilon_{i-1}x)((\sigma_0\dots \sigma_{k}\sigma_{k+1}\dots\sigma_{r}))\ \hbox{ by definition of $\epsilon_i$} &&\\
&= (h^-(\epsilon_{i-1}x))(k_1\dots k_{n+2})\ \hbox{ by definition of $h^-$} &&
\end{alignat*}
\item $k_i\neq  k_{i+1}$. Now $s=r$ and since $\{k_i,k_{i+1}\}\subset \{-,0\}$,
then there exists a unique $j_0$ such that $\sigma_{j_0}+1\in\{i,i+1\}$ and we have
\beas
&&\sigma_j+1\leq i-1 \Longrightarrow \sigma_j=\tau_j \hbox{ (in this case, $j<j_0$)}\\
&& \tau_{j_0}+1=i \\
&& \sigma_j+1\geq i+2 \Longrightarrow \sigma_j- 1=\tau_{j} \hbox{ (in this case, $j>j_0$)}
\eeas
There are two subcases : $\sigma_{j_0}+1=i$ and $\sigma_{j_0}+1=i+1$. In the
first situation,

{
\begin{alignat*}{2}
& x((\tau_0\dots \tau_s))&&\ \\
&=x((\sigma_0<\dots <\sigma_{j_0-1}<\sigma_{j_0}=i-1<\sigma_{j_0+1}-1<\dots<\sigma_r-1))&&\\
&=x((\sigma_0<\dots <\sigma_{j_0-1}<\sigma_{j_0}<\sigma_{j_0+1}-1<\dots<\sigma_r-1))&&\ \\
&= (\epsilon_{i-1}x)((\sigma_0<\dots <\sigma_{j_0}<\sigma_{j_0+1}<\dots<\sigma_r))\ \hbox{ by definition of $\epsilon_i$} &&\\
&= (h^-(\epsilon_{i-1}x))(k_1\dots k_{n+2})\ \hbox{ by definition of $h^-$}&&
\end{alignat*}}

In the second situation,

{
\begin{alignat*}{2}
& x((\tau_0\dots \tau_s))\\
&=x((\sigma_0<\dots <\sigma_{j_0-1}<\sigma_{j_0}-1=i-1<\sigma_{j_0+1}-1<\dots<\sigma_r-1)) &&\\
&=x((\sigma_0<\dots <\sigma_{j_0-1}<\sigma_{j_0}-1<\sigma_{j_0+1}-1<\dots<\sigma_r-1))&&\ \\
&= (\epsilon_{i-1}x)((\sigma_0<\dots <\sigma_{j_0}<\sigma_{j_0+1}<\dots<\sigma_r))\ \hbox{ by definition of $\epsilon_i$} &&\\
&= (h^-(\epsilon_{i-1}x))(k_1\dots k_{n+2})\ \hbox{ by definition of $h^-$}&&
\end{alignat*}}

\end{enumerate}

\end{proof}

Notice that $h^-$ induces a natural transformation from $\CR^{gl}_*$
to $\CR^{-}_*$ which is not injective. Consider for example the
$\omega$-category consisting of two composable $1$-morphisms $u$ and
$v$ with $t_0u=s_0v$. The $0$-simplexes $u$ and $u*_0v$ of
$\mathcal{N}^{gl}_0$ have indeed the same image by $h^-$ in
$\CR^{-}_1$. To see that, consider the thin square $c$  from $I^2$ to
$\C$ defined by $c(-0)=u*_0v$, $c(0+)=t_0v$, $c(0-)=u$, $c(+0)=v$ and
$c(00)=u*_0v$.

Now we arrive at :

\bth There exists one and only one morphism of cuts  from
$\mathcal{N}^{gl}$ to $\mathcal{N}^{-}$. \eth

The proof of this theorem uses Theorem~\ref{grading}
assertion~\ref{extremite} as shortcut. There is no vicious circle
because the uniqueness of $h^-$ and $h^+$ is used nowhere in this
paper. The only fact which is used is that Theorem~\ref{simplicial_h}
provides a natural transformation from $\mathcal{N}^{gl}$ to
$\mathcal{N}^{-}$ which is injective on the underlying sets.

\bpf Let $h$ and $h'$ be two morphisms of cuts from
$\mathcal{N}^{gl}$ to $\mathcal{N}^{-}$. One proves by induction on $n$
that $h_n$ and $h'_n$ from $\mathcal{N}^{gl}_n$ to $\mathcal{N}^{-}_n$
coincide. For $n=0$, $\mathcal{N}^{gl}_0=\mathcal{N}^{-}_n=tr^0\P$. The only
natural transformation from $tr^0\P$ to itself is $Id_{tr^0\P}$, therefore
$h_0=h'_0$.

Suppose $P(n)$ proved for
some $n\geq 0$. Then for any $x\in \mathcal{N}^{gl}_{n+1}(\C)$, and
for any $0\leq i \leq n+1$,
\begin{alignat*}{2}
 \de_{i+1}^-h_{n+1}(x) &= h_{n}(\de_i x) &&\ \hbox{ since $h$ morphism of simplicial sets}\\
&= h'_{n}(\de_i x) &&\ \hbox{ by induction hypothesis}\\
&= \de_{i+1}^- h'_{n+1}(x) &&\ \hbox{ since $h'$ morphism of simplicial sets}
\end{alignat*}
Now with $1\leq j\leq n+2$,
\begin{alignat*}{2}
&(\de_{j}^+h_{n+1}(x))(-_{n+1})&&\\
&=h_{n+1}(x)(-\dots -[+]_j -\dots -)&&\\
&= h_{n+1}(x)\left( t_0 R(-\dots -[0]_j -\dots -)\right) &&\\
&= t_0 \left(h_{n+1}(x)(R(-\dots -[0]_j -\dots -))\right)\ \hbox{ since $h_{n+1}(x)$ $\omega$-functor}&&\\
&= t_0 \left( (\de_1^-\dots \widehat{\de_j^-}\dots \de_{n+2}^- h_{n+1}(x))(0) \right)&&\\
&= t_0 \left( h_0(\de_0\dots \widehat{\de_{j-1}}\dots \de_{n+1}x) (0)\right)\ \hbox{ since $h$ morphism of simplicial sets}&&\\
&= t_0 \left((\de_0\dots \widehat{\de_{j-1}}\dots \de_{n+1}x)((0))\right)&&
\end{alignat*}
So the $0$-morphism $\de_{j}^+h_{n+1}(x))(-_{n+1})$ is the value of
the constant map $t_0\circ x$ of Theorem~\ref{grading} (denoted by
$T(x)$ in Section~\ref{deformation}).

Let $\D$ be the unique $\omega$-category
such that $\P\D=\Delta^{n+1}$ and with $\D_0=\{\alpha,\beta\}$,
$s_0(\P\D)=\{\alpha\}$, $t_0(\P\D)=\{\beta\}$ and $\alpha\neq \beta$. And consider
$Id_{\Delta^{n+1}}\in \mathcal{N}^{gl}_{n+1}(\D)$.

Suppose that $+\in\{k_1,\dots ,k_{n+2}\}\subset \{-,+\}$ and suppose
that at least two $k_i$ are equal to $+$. Then there exists a
$1$-morphism $u$ of $I^{n+2}$ such that $s_0u=\ell_1\dots \ell_{n+2}$ with
exactly one $\ell_i$ equal to $+$ and such that $t_0u=k_1\dots k_{n+2}$.
Then
\[s_0\left(h_{n+1}(Id_{\Delta^{n+1}})(u)\right)= h_{n+1}(Id_{\Delta^{n+1}})(\ell_1\dots \ell_{n+2})=\beta\]
by the previous calculation.
Since $\beta$ is the unique morphism of $\D$ with  $0$-source $\beta$, then
$h_{n+1}(Id_{\Delta^{n+1}})(u)=\beta$ and therefore
$$h_{n+1}(Id_{\Delta^{n+1}})(k_1\dots k_{n+2})=\beta.$$

Suppose now that $+\in\{k_1,\dots ,k_{n+2}\}$ with perhaps  some $0$ in the set. Then
\[s_0\left(h_{n+1}(Id_{\Delta^{n+1}})(k_1\dots k_{n+2})\right)=\beta\] and
therefore
\[\ev \circ h_{n+1}(Id_{\Delta^{n+1}})(k_1\dots k_{n+2})=\beta=T\left(Id_{\Delta^{n+1}}\right).\]
The $\omega$-functor $x$ from
$\Delta^{n+1}$ to $\P\C$ induces a non-contracting $\omega$-functor $\overline{x}$
from  $\D$ to $\C$ with $\overline{x}(\alpha)=S(x)$ ($S(x)$ being the value
of the constant map $s_0\circ x$ by Theorem~\ref{grading}) and
$\overline{x}(\beta)=T(x)$
which sends $Id_{\Delta^{n+1}}\in \mathcal{N}^{gl}_{n+1}(\D)$
on $x\in \mathcal{N}^{gl}_{n+1}(\C)$. So by naturality,
\[\ev \circ h_{n+1}(x)(k_1\dots k_{n+2})=T(x).\]
Therefore for any $1\leq j\leq n+2$,
$\de_{j}^+h_{n+1}(x)=\de_{j}^+h'_{n+1}(x)$. By hypothesis,
$\ev(h_{n+1}(x))=\ev(x)=\ev(h'_{n+1}(x))$. So $h_{n+1}(x)$ and $h'_{n+1}(x)$
induce the same labeling of the faces of $I^{n+2}$ and $P(n+1)$ is proved.
\epf

Without explanation, here is the construction of $h^+$ :

\bp\label{simplicial_h_plus} Let $x$ be an $n$-simplex of the globular simplicial nerve of
$\C$. Then the map $h_n^+(x)$ from $\Cube^{n+1}$ to $\C$ defined by
\begin{enumerate}
\item $-\in \{k_1\dots k_{n+1}\}$ implies
$h_n^+(x)(k_1\dots k_{n+1})=s_0 x((n))$ (notice that $(n)$ is the initial state of
$\Delta^n$)
\item $\{k_1,\dots ,k_{n+1}\}\subset \{+,0\}$ and $$\{k_1,\dots, k_{n+1}\}\cap \{0\}=
\{k_{\sigma_0+1},\dots,k_{\sigma_r+1}\}$$ with $\sigma_0<\dots<\sigma_r$ implies
$h_n^+(x)(k_1\dots k_{n+1})=x((\sigma_0\dots\sigma_r))$
\item $h_n^+(x)(+_{n+1})=t_0 x((0))$ (notice that $(0)$ is the final state of
$\Delta^n$)
\end{enumerate}
yields an $\omega$-functor from $I^{n+1}$ to $\C$. Moreover, $h^+$
induces a morphism of simplicial sets from the globular nerve of $\C$
to its positive corner nerve. And the map from
$\mathcal{N}^{gl}_{-1}(\C)$ to $\mathcal{N}^{+}_{-1}(\C)$ defined by
$(x,y)\mapsto y$ extends the previous morphism to the corresponding
augmented simplicial nerves. Moreover for $n\geq 0$, $h_n^+$ is a
one-to-one map and the image of $h_n^+$ contains exactly all cubes $x$
of the positive corner nerve such that as soon as $\de_i^-x$ exists,
then it is $0$-dimensional. \ep

\begin{question} Is it possible to find an appropriate setting
where the globular cut would be an initial object ? Is it possible to
characterize the diagram of cuts of Figure~\ref{fundamental} ?
\end{question}

As immediate corollary of the construction of $h^-$ and its
injectivity, let us introduce the analogue of
Proposition~\ref{filling} in the globular nerve.

\bd In a simplicial set $A$, a $n$-shell is a family $(x_i)_{i=0,\dots,n+1}$
of $(n+2)$ $n$-simplexes of $A$ such that for any $0\leq i<j\leq n+1$,
$\de_i x_j=\de_{j-1} x_i$. \ed

\bp\label{filling_simp} Let $\C$ be a non-contracting $\omega$-category. Consider
a $n$-shell $(x_i)_{i=0,\dots,n+1}$ of the globular simplicial nerve
of $\C$. Then
\begin{enumerate}
\item The labeling defined by $(x_i)_{i=0,\dots,n+1}$ yields an
$\omega$-functor $x$ (and necessarily exactly one)
from $\Delta^{n+1}\backslash \{(01\dots n+1)\}$ to $\P\C$.
\item Let $u$ be a morphism of $\C$ such that
\[s_{n}u=x\left(s_n R((01\dots n+1))\right)\]
and
\[t_{n}u=x\left(t_n R((01\dots n+1))\right)\]
Then there exists one and only
one $\omega$-functor still denoted by $x$ from $\Delta^{n+1}$ to $\P\C$ such that
for any $0\leq i\leq n+1$, $\de_i x = x_i$ and $$x((01\dots n+1))=u.$$
\end{enumerate}
\ep

\section{Regularity of the globular cut}\label{regularity}

This section is devoted to the proof of the following theorem.

\bth\label{normal} The globular cut is regular. \eth

The principle of this proof is to use the injectivity of the natural
transformation $h^-$ from $\mathcal{N}^{gl}$ to $\mathcal{N}^{-}$ and
to use the regularity of $\mathcal{N}^{-}$.

The folding operator $\Phi_n^{gl}:=\Phi_n^{\mathcal{N}^{gl}}$ is
called the $n$-dimensional globular folding operator and we set
$\square_n^{gl}:=\square_n^{\mathcal{N}^{gl}}$. It is clear that
rule~\ref{regular0} and rule~\ref{regular1} of Definition~\ref{def_regular}
are satisfied. We have to check the rest of it.

\bth
For any natural transformation of functors $\mu $ from $\mathcal{N}^{gl}_{n-1}$ to
 $\mathcal{N}^{gl}_{n}$ with $n\geq 1$, and for any natural map $\square$ from
$tr^{n-1}\P$ to $\mathcal{N}^{gl}_{n-1}$ such that
$\ev\circ \square=\id_{tr^{n-1}\P}$, there exists one and only one
natural transformation denoted by $\mu .\square$ from $tr^n\P$ to $\mathcal{N}^{gl}_{n}$ such
that the following diagram commutes
\[\xymatrix{\ar@/^20pt/[rr]^{\id_{tr^n\P}}{tr^n \P}\fr{\mu .\square}& {\mathcal{N}^{gl}_n}\fr{\ev}&{tr^n \P}\\
\ar@/_20pt/[rr]_{\id_{tr^{n-1}\P}}{tr^{n-1}\P}\fu{i_{n}}\fr{\square}&{\mathcal{N}^{gl}_{n-1}}\fr{\ev}\fu\mu &{tr^{n-1}\P}\fu{i_{n}}}\]
where $i_{n}$ is the canonical inclusion functor from
$tr^{n-1}\P$ to $tr^n \P$. \eth

\bpf The natural transformation $h^- \square$ from $tr^{n-1}\P$ to
$\mathcal{N}^-_{n-1}$ can be lifted to a natural transformation
$(h^-(\mu)).(h^-\square)$ from $tr^{n}\P$ to $\mathcal{N}^-_{n}$
since the cut $\mathcal{N}^-$ is regular. Since
$h^-(\mu.\square)=(h^-(\mu)).(h^-\square)$ and since $h^-$ is one-to-one
in positive degree, there is at most one solution for  this lifting
problem.

\[\xymatrix{\ar@/^20pt/[rr]^{h^-(\mu).(h^-\square)}{tr^n \P}& {\mathcal{N}^{gl}_n}\fr{h^-}&{\mathcal{N}^{-}_n}\\
\ar@/_20pt/[rr]_{h^- \square}{tr^{n-1}\P}\fu{i_{n}}\fr{\square}&{\mathcal{N}^{gl}_{n-1}}\fr{h^-}\fu{\mu} &{\mathcal{N}^{-}_{n-1}}\fu{h^-(\mu)}}\]
Let $x\in\C_{n+1}$. For
$0\leq i\leq n$, the natural transformation
$$\ev\ \de_i\left(h^-(\mu).(h^-\square)\right):tr^n \P \rightarrow tr^{n-1}\P$$
is of the
form $d_{m_i}^{\alpha_i}$ for some $\alpha_i\in\{-,+\}$ and some $m_i\leq n$.
Therefore
\begin{alignat*}{2}
&\de_i\left(h^-(\mu).(h^-\square)\right)&&\\
&=\de_i\left(h^-(\mu).(h^-\square)\right)i_n d_{m_i}^{\alpha_i}&&\ \hbox{ by Definition~\ref{def_regular} rule~\ref{regular3.02}}\\\
&=\de_ih^-(\mu)h^-\square d_{m_i}^{\alpha_i}&&\ \hbox{ by hypothesis}\\
&=\de_ih^-\mu\square d_{m_i}^{\alpha_i}&&\\
&=h^- \de_i\mu\square d_{m_i}^{\alpha_i}&&\ \hbox{ since $h^-$ morphism of simplicial sets}
\end{alignat*}
So $\de_i\left(h^-(\mu).(h^-\square)\right)(x)\in h^-(\mathcal{N}_{n-1}^{gl}(\C))$
for any $0\leq i\leq n$ and by Proposition~\ref{filling_simp},
$\left(h^-(\mu).(h^-\square)\right)(x)\in h^-(\mathcal{N}_n^{gl}(\C))$. Let
$\square '(x)$ be the unique element of $\mathcal{N}_n^{gl}(\C)$ such that
\[h^- \square '(x):= \left(h^-(\mu).(h^-\square)\right)(x)\]
Then $\square '$ is a solution.
\epf

\begin{cor} The equalities $h^- \Phi^{gl}=\Phi^- h^-$ and
$h^+ \Phi^{gl}=\Phi^+ h^+$ hold. \end{cor}

\bpf It is a consequence of the naturality of $h^-$ and $h^+$
and of Proposition~\ref{Phi_functoriel}. \epf

Now here is  a characterization of globular folding operators :

\bp Let $x$ be a $n$-simplex of the globular nerve of $\C$. Then
$x=\Phi^{gl}(x)$ if and only if for $0\leq i\leq n-2$,
$\de_i x\in Im(\epsilon_{n-2}\dots \epsilon_i)$. \ep

\bpf The equality $x=\Phi^{gl}(x)$ implies
$h^-(x)=\Phi^{-}(h^-(x))$, implies by Theorem~\ref{caracterisation} that
for $1\leq i\leq n-1$,
\beas
h^-(\de_{i-1}x)&=&\de^-_i(h^-(x))= \Gamma_{n-1}^-\dots \Gamma_{i}^- \square_i^- d_i^{(-)}h^-(x)(0_{n+1})\\
&=&h^-\left(\epsilon_{n-2}\dots \epsilon_{i-1} \square_i^{gl}s_i x((0\dots n))\right)
\eeas
therefore $\de_{i-1}x\in Im(\epsilon_{n-2}\dots
\epsilon_{i-1})$. Conversely, if for $0\leq i\leq n-2$,
$\de_i x\in Im(\epsilon_{n-2}\dots \epsilon_i)$, then
$h^-(x)=\Phi^- h^-(x)=h^-\Phi^{gl}(x)$ and therefore
$x=\Phi^{gl}(x)$.
\epf

\bth The globular folding operator $\Phi^{gl}$ induces the identity
map on the globular reduced chain complex $\CR_*^{gl}$. \eth

\bpf Consider the $\theta_i^-$ operators of Theorem~\ref{decomposition_theta}. If $x\in
\mathcal{N}_n^{gl}$, then $h^-x$ is negative. So $\theta_i^-h^-x$ is
also negative by Theorem~\ref{decomposition_theta}(\ref{th1}) and
determines a unique element $\theta_i^{gl}x\in \mathcal{N}_n^{gl}$
such that $h^-\theta_i^{gl}x=\theta_i^-h^-x$. It is clear that these
operators $\theta_i^{gl}$ induces the identity map on the reduced
globular complex by Theorem~\ref{decomposition_theta}(\ref{th2}).  Since
$\Phi^-h^-x$ is also negative, then by
Theorem~\ref{decomposition_theta}(\ref{th3}),
\[\Phi^-h^-x=\theta_{i_1}^- \dots \theta_{i_s}^- h^-x\]
for some sequence $i_1,\dots,i_s$. Therefore by the injectivity of $h^-$,
\[\Phi^{gl}x=\theta_{i_1}^{gl} \dots \theta_{i_s}^{gl} x\]
\epf

\bth\label{relation_glob} In the reduced globular complex, one has
\[\square_n^{gl}(x*_p y)=\square_n^{gl}(x)+ \square_n^{gl}(y)\]
for any morphisms $x$ and $y$ of $\C$ of dimension $n$ and for
$1\leq p\leq n-1$.
\eth

\begin{proof}[Sketch of proof] One has
\beas
h^-(\square_n^{gl}(x*_p y))&=&\square_{n}^- (x *_p y)\\
&=& \square_{n}^- (x)+ \square_{n}^- (y)+t_1 +\de^- t_2\\
&=&h^-(\square_n^{gl}(x))+h^-(\square_n^{gl}(y))+t_1 +\de^- t_2
\eeas
 with $t_1$ a thin $(n+1)$-cube and $t_2$ a thin $(n+2)$-cube. The
 proof made in \cite{Coin} shows that $t_1$ and $t_2$ are in the image
 of $h^-$. Indeed, the existence of $t_1$ and $t_2$ comes from the
 vanishing of some globular nerve. Therefore $t_1=h^-(T_1)$ and
 $t_2=h^-(T_2)$ where $T_1$ is a thin $n$-simplex and $T_2$ a thin
 $(n+1)$-simplex. This completes the proof.
\end{proof}

In fact one can explicitly verify that if $x$ and $y$ are two
$n$-morphisms of $\C$, then $\square_{n}^{gl}(x *_{n-1}
y)-\square_{n}^{gl}(x)-\square_{n}^{gl}(y)$ is a boundary in the
normalized globular complex. It suffices to consider the thin
$(n+1)$-cube $B^n_{n-1}(x,y)$ of \cite{Coin} which turns to be in
the image of $h^-$ because it is negative. Therefore with
$b(x,y)\in \omega Cat(\Delta^n,\P\C)$ defined by $\de_i
b(x,y)=\epsilon_{n-2}\dots \epsilon_i \square_{i+1}^{gl}
d_{i+1}^{(-)^{i+1}}x$ for $0\leq i\leq n-3$ (observe that
$d_{i+1}^{(-)^{i+1}}x= d_{i+1}^{(-)^{i+1}}y$),
$\de_{n-2}b(x,y)=\square_{n}^{gl}y$,
$\de_{n-1}b(x,y)=\square_{n}^{gl}(x*_{n-1} y)$, $\de_n
b(x,y)=\square_{n}^{gl}x$, one has
\[\de b(x,y)=\pm \left(\square_{n}^{gl}(x *_{n-1}
y)-\square_{n}^{gl}(x)-\square_{n}^{gl}(y)\right)+\hbox{ degenerate elements}.\]

\section{Example of calculations of globular homology}\label{vanish_In}

The main goal of this section is to prove the vanishing of the
globular homology of the $n$-cube in positive dimension for all $n\geq
0$. However we also study the case of the $\omega$-category $2_n$ generated
by one $n$-morphism and pose some questions about the globular homology
of the $\omega$-category generated by a composable pasting scheme in the
sense of \cite{CPS}.

\bth For any $p>0$ and any $n\geq 0$, $H_p^{gl}(2_n)=0$. \eth

\bpf For $p=1$, it is obvious. For $p>1$, one has
$$H_p^{gl}(2_n)\iso H_{p-1}(\P 2_n) \iso H_{p-1}(2_{n-1})=0$$
where $H_*(\D)$ means the
simplicial homology of the simplicial nerve of the $\omega$-category
$\D$. \epf

\bd \cite{Gau} Let $\C$ be an $\omega$-category and let $\alpha$ and
$\beta$ be two $0$-morphisms of $\C$. Then the \textit{bilocalization}
of $\C$ with respect to $\alpha$ and $\beta$ is the
$\omega$-subcategory of $\C$ obtained by keeping in dimension $0$ only
$\alpha$ and $\beta$ and by keeping in positive dimension all morphisms
$x$ such that $s_0x=\alpha$ and $t_0x=\beta$. It is denoted by
$\C[\alpha,\beta]$. \ed

\bth\label{grading} Let $\C$ be a non-contracting  $\omega$-category.
\begin{enumerate}
\item\label{extremite} Let $x$ be an $\omega$-functor from $\Delta^n$ to $\P\C$ for some
$n\geq 0$. Then the set maps $$(\sigma_0\dots \sigma_r)\mapsto s_0 x((\sigma_0\dots \sigma_r))$$
and $$(\sigma_0\dots \sigma_r)\mapsto t_0 x((\sigma_0\dots \sigma_r))$$ from
the underlying set of faces of $\Delta^n$ to $\C_0$ are constant. The unique value of
$s_0\circ x$ is denoted by $S(x)$ and the unique value of $t_0 \circ x$ is denoted by
$T(x)$.
\item For any pair $(\alpha,\beta)$ of $0$-morphisms of $\C$, for any
$n\geq 1$, and for any $0\leq i\leq n$, then
$\de_i \left(\mathcal{N}_n^{gl}(\C[\alpha,\beta])\right)\subset \mathcal{N}_{n-1}^{gl}(\C[\alpha,\beta])$.
\item For any pair $(\alpha,\beta)$ of $0$-morphisms of $\C$, for any
$n\geq 0$, and for any $0\leq i\leq n$, then
$\epsilon_i\left(\mathcal{N}_n^{gl}(\C[\alpha,\beta])\right) \subset \mathcal{N}_{n+1}^{gl}(\C[\alpha,\beta])$.
\item By setting, $G^{\alpha,\beta}\mathcal{N}_n^{gl}(\C):=\mathcal{N}_n^{gl}(\C[\alpha,\beta])$
for $n\geq 0$ and $G^{\alpha,\beta}\mathcal{N}_{-1}^{gl}(\C)\linebreak[4]:=\{(\alpha,\beta),(\beta,\alpha)\}$,
one
obtains a $(\C_0\p \C_0)$-graduation on the globular nerve ; in particular, one has the direct sum
of augmented simplicial sets
\[\mathcal{N}_*^{gl}(\C)=\bigsqcup_{(\alpha,\beta)\in\C_0\p\C_0} G^{\alpha,\beta}\mathcal{N}_*^{gl}(\C)\]
and $G^{\alpha,\beta}\mathcal{N}_*^{gl}(\C)=\mathcal{N}_*^{gl}(\C[\alpha,\beta])$.
\end{enumerate}
\eth

\bpf The only non-trivial part is the first assertion. Let $P(n)$ be the
property : ``for any non-contracting $\omega$-category $\C$ and any
$\omega$-functor $x$ from $\Delta^n$ to $\P\C$, the set map
$(\sigma_0\dots \sigma_r)\mapsto s_0 x((\sigma_0\dots \sigma_r))$ from
the set of faces of $\Delta^n$ to $\C_0$ is constant.''

There is nothing to check for $P(0)$. For $P(1)$, if $x$ is an
$\omega$-functor from $\Delta^1$ to $\P\C$, then $s_1 x((01))=x((1))$
and $t_1 x((01))=x((0))$ in $\C$. Therefore $$s_0 x((01)) = s_0 s_1
x((01))=s_0 x((1))$$ and $$s_0 x((0))=s_0 t_1 x((01))= s_0 x((01)).$$
Therefore $P(1)$ is true.

Suppose $P(n)$ proved for some $n\geq 1$ and let us prove $P(n+1)$.
For any $1\leq i \leq n$, the $\omega$-functor
$x:\Delta^{n+1}\rightarrow \P\C$ induces an $\omega$-functor on the
$\omega$-category $\Delta^{n+1}_i$ generated by the face
$(0\dots\widehat{i}\dots n+1)$ and its subfaces. One has an
isomorphism of $\omega$-categories $\Delta^n\iso \Delta^{n+1}_i$.
Therefore the restriction of $s_0\circ x$ to the faces of
$\Delta^{n+1}_i$ is constant by induction hypothesis. Now it is clear
that $\Delta^{n+1}_i\cap \Delta^{n+1}_{i+1}\iso \Delta^{n-1}\neq
\emptyset$ since $n\geq 1$. Therefore the set map $s_0\circ x$
restricted to $\Delta^{n+1}_i\cup \Delta^{n+1}_{i+1}$ is constant.
Therefore the restriction of the set map $s_0\circ x$ to the faces of dimension at most
$n$ of $\Delta^{n+1}$ is constant. We know that
\[s_n R((01\dots n+1))=\Psi(X_0,X_1,\dots,X_s)\]
where $X_0,X_1,\dots,X_s$ are faces of $\Delta^{n+1}$ of dimension at
most $n$. So
\begin{alignat*}{2}
s_0 x((01\dots n+1))&= s_0 s_{n+1} x ((01\dots n+1))&&\\
&= s_0 x \left(s_n R((01\dots n+1))\right) &&\ \hbox{ since $x$ $\omega$-functor}\\
&= s_0 x \Psi(X_0,X_1,\dots,X_s)
\end{alignat*}
where $\Psi$ is a function using only the compositions of $\Delta^{n+1}$. Then
\[x \Psi(X_0,X_1,\dots,X_s)= \Psi'(x(X_0), x(X_2),\dots, x(X_s))\]
where $\Psi'$ is obtained
from $\Psi$ by replacing $*_i$ by $*_{i+1}$ since $x$ is an
$\omega$-functor from $\Delta^{n+1}$ to $\P\C$. So
\[s_0 x((01\dots n+1))=  \Psi'(s_0 x(X_0), s_0 x(X_2),\dots, s_0 x(X_s))=s_0 x(X_0) \]
with the axioms of $\omega$-categories. Therefore $P(n+1)$ is proved.
\epf

\bd Let $\C$ be a non-contracting $\omega$-category with exactly one
initial state $\alpha$ and one final state $\beta$. Then the bilocalization
$\C[\alpha,\beta]$ is also non-contracting and one can set
$\Omega\C= \P (\C[\alpha,\beta])$. \ed

\bth\cite{MR34:6767,MR81m:55010,MR91m:32042}
\label{derive} Let $n\geq 1$. Then $\Omega \Delta^n=I^{n-1}$
and $\Omega I^{n-1}=P^{n-1}$ where $P^{n-1}$ is the free
$\omega$-category generated by the composable pasting scheme of the
faces of the $(n-1)$-dimensional permutohedron. \eth

\bth\label{cube} For any $n\geq 0$, and any $p>0$, $H_p^{gl}(I^n)=0$. \eth

\bpf One has $H_p^{gl}(I^n)=\bigoplus_{(\alpha,\beta)\in
  \C_0\p\C_0}H_p^{gl}(I^n[\alpha,\beta])$ by Theorem~\ref{grading}. So
it suffices to prove the vanishing of $H_p^{gl}(I^n[\alpha,\beta])$ as
soon as $I^n[\alpha,\beta]$ contains morphisms in strictly positive
dimension to prove the theorem.

Let $\alpha$ and $\beta$ be two $0$-morphisms of $I^n$ such that
$I^n[\alpha,\beta]$ contains other morphisms than $\alpha$ and
$\beta$.  Then in particular it contains some $1$-morphisms from
$\alpha$ to $\beta$ which is a composite of $1$-dimensional faces of
$I^n$. Suppose that $\alpha=k_1\dots k_n$. Then $\beta$ is obtained
from $\alpha$ by replacing some $k_i$ equal to $-$ by $+$. Let
$k_{\sigma_1},\dots ,k_{\sigma_r}$ be these $k_i$. Then
$$I^n[\alpha,\beta]\iso I^r[-_r,+_r]$$ as $\omega$-category. Therefore
it suffices to prove that $H_p^{gl}(I^n[-_n,+_n])$ vanishes.

The vanishing of $H_1^{gl}(I^n[-_n,+_n])$ is obvious.
One has $$H_p^{gl}(I^n[-_n,+_n])=H_{p-1}(P^n)$$ for $p\geq 2$ by
Theorem~\ref{derive} and $H_{p-1}(P^n)=0$ because the simplicial nerve
of a composable pasting scheme is contractible \cite{CPS}.
\epf

\bth\label{delta} For any $n\geq 0$, and any $p>0$, $H_p^{gl}(\Delta^n)=0$. \eth

\bpf By proceeding as in Theorem~\ref{cube}, we see that it suffices
to prove that \[H_p^{gl}(\Delta^n[(r),(s)])=0\] for any pair $((r),(s))$
of $0$-morphisms of $\Delta^n$ and for $n\geq 2$. However,
$\Delta^n[(r),(s)]$ is non-empty if and only if $r>s$ with our
conventions and in this case,
$$\Delta^n[(r),(s)]\iso\Delta^{r-s}[(r-s),(0)].$$ Therefore
$H_p^{gl}(\Delta^n[(r),(s)])\iso H_{p-1}(I^{r-s-1})$ by Theorem~\ref{derive}.
\epf

More generally, as in \cite{Coin}, one sees that if $\C$ is a
non-contracting $\omega$-category such that $\P\C$ is the free
$\omega$-category generated by a composable pasting scheme in the
sense of \cite{CPS}, then $H_p^{gl}(\C)=0$ for $p\geq 1$. This is
related to the problem of the existence of the derived pasting scheme
of a given composable pasting scheme \cite{MR92j:18004}.

\begin{conj} Let $\C$ be an $\omega$-category which is the free $\omega$-category
  generated by a composable pasting scheme (therefore $\C$ is
  non-contrac\-ting).  Then for any $p>0$, $H_p^{gl}(\C)=0$.
\end{conj}

\section{Relation between the new globular homology and the old one}\label{relation}

First of all, recall the definition of both formal corner homology
theories from \cite{Coin}.

\bd Let $\C$ be a non-contracting $\omega$-category. Set
\begin{itemize}
\item $\CF^{-}_0(\C):=\Z\C_0$
\item $\CF^{-}_1(\C):=\Z\C_1$
\item $\CF^{-}_n(\C)= \Z\C_n/\{x*_0y=x,x*_1y=x+y, \dots , x
*_{n-1}y=x+y\hbox{ mod }\Z tr^{n-1}\C\}$ for $n\geq 2$
\end{itemize}
with the differential map $s_{n-1}-t_{n-1}$ from $\CF^{-}_n(\C)$ to
$\CF^{-}_{n-1}(\C)$ for $n\geq 2$ and $s_0$ from
$\CF^{-}_1(\C)$ to $\CF^{-}_0(\C)$. This chain complex is called
the formal negative corner complex. The associated homology is denoted by
$\HF^{-}(\C)$ and is called the formal negative corner homology of $\C$.
The map $\CF_*^-$ (resp. $\HF_*^-$) induces a functor from $\omega Cat_1$
to $\comp$ (resp. $Ab$).
\ed

and symmetrically

\bd Let $\C$ be a non-contracting $\omega$-category. Set
\begin{itemize}
\item $\CF^{+}_0(\C):=\Z\C_0$
\item $\CF^{+}_1(\C):=\Z\C_1$
\item $\CF^{+}_n(\C)= \Z\C_n/\{x*_0y=y,x*_1y=x+y, \dots , x
*_{n-1}y=x+y\hbox{ mod }\Z tr^{n-1}\C\}$ for $n\geq 2$
\end{itemize}
with the differential map $s_{n-1}-t_{n-1}$ from $\CF^{+}_n(\C)$ to
$\CF^{+}_{n-1}(\C)$ for $n\geq 2$ and $t_0$ from
$\CF^{+}_1(\C)$ to $\CF^{+}_0(\C)$. This chain complex is called
the formal positive corner complex. The associated homology is denoted by
$\HF^{+}(\C)$ and is called the formal positive corner homology of $\C$.
The map $\CF_*^+$ (resp. $\HF_*^+$) induces a functor from $\omega Cat_1$
to $\comp$ (resp. $Ab$).
\ed

The maps $\square_n^\pm$ from $\C_n$ to $C_n^\pm(\C)$ induce a natural
transformation from $\CF_*^\pm$ to $\CR_*^\pm$ and a natural transformation
from $\HF_*^\pm$ to $\HR_*^\pm$.

\bd\label{new_glob} Let $\C$ be a non-contracting $\omega$-category. Set
\begin{itemize}
\item $\CF^{gl}_0(\C):=\Z\C_0\ot\Z\C_0\iso \Z (\C_0\p\C_0)$
\item $\CF^{gl}_1(\C):=\Z\C_1$
\item $\CF^{gl}_n(\C)= \Z\C_n/\{x*_1y=x+y, \dots , x
*_{n-1}y=x+y\hbox{ mod }\Z tr^{n-1}\C\}$ for $n\geq 2$
\end{itemize}
with the differential map $s_{n-1}-t_{n-1}$ from $\CF^{gl}_n(\C)$ to
$\CF^{gl}_{n-1}(\C)$ for $n\geq 2$ and $s_0\ot t_0$ from
$\CF^{gl}_1(\C)$ to $\CF^{gl}_0(\C)$. This chain complex is called
the formal globular complex. The associated homology is denoted by
$\HF^{gl}(\C)$ and is called the formal globular homology of $\C$.
\ed

By Theorem~\ref{relation_glob} and Corollary~\ref{diff_reduite}, we
see that the globular folding operators induce a natural morphism of
chain complex from $\CF_*^{gl}$ to $\CR_*^{gl}$, and therefore a
natural transformation from $\HF_*^{gl}$ to $\HR_*^{gl}$.

\begin{question} When is the natural morphism of chain complexes
$R^{gl}$ from $\CF_*^{gl}(\C)$ to $\CR_*^{gl}(\C)$ a quasi-isomorphism ?
\end{question}

\begin{conj}\label{thin_glob} (About the thin elements of the globular
complex) Let $\C$ be a globular $\omega$-category which is either the
free globular $\omega$-category generated by a semi-cubical set or the free
globular $\omega$-category generated by a globular set.  Let $x_i$ be
elements of $C_{n}^{gl}(\C)$ and let $\lambda_i$ be natural
numbers, where $i$ runs over some set $I$. Suppose that for any $i$,
$\ev(x_i)$ is of dimension strictly lower than $n$ (one calls it a
thin element). Then $\sum_i \lambda_i x_i$ is a boundary if and only
if it is a cycle.  \end{conj}

The above conjecture is clear for $C_{2}^{gl}$ because all thin elements 
are degenerate. In higher dimension, there is enough room to have thin 
elements which are composition of degenerate elements, but which are not 
degenerate themselves.

The above conjecture is equivalent to claiming that the globular
homology and the reduced one are equivalent for free globular
$\omega$-categories generated by either a semi-cubical set or a globular
set.

Now we are in position to give the exact statement relating
the old globular homology of \cite{Gau} and the new one.

\bd\cite{Gau} Let $(C_*^{old-gl}(\C),\de^{old-gl})$ be the chain
complex defined as follows : \linebreak[4] $C_0^{old-gl}(\C)=\Z\C_0\oplus \Z\C_0$ and
for $n\geqslant 1$, $C_n^{old-gl}(\C)=\Z\C_n$, $\de^{old-gl}(x)=(s_0x,t_0x)$
if $x\in\Z\C_1$ and for $n\geqslant 1$, $x\in\Z\C_{n+1}$ implies
$\de^{old-gl}(x)=s_nx-t_nx$. This complex is called the \textit{old globular
  complex} of $\C$ and its corresponding homology the old globular
homology.  \ed

\begin{figure}
\[
\xymatrix{ && \\
\ar@{->}[ru]^{u}\ar@{->}[rd]^{w} && \ar@{->}[lu]^{v}\ar@{->}[ld]^{x}\\
&&
}
\]
\caption{A false $1$-globular cycle in the old globular homology}
\label{false}
\end{figure}
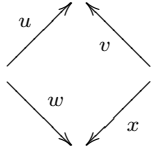

Instead of $C_0^{old-gl}(\C)=\Z\C_0\oplus \Z\C_0$, we set
$C_0^{old-gl}(\C)=\Z(\C_0\otimes\C_0)$ with the differential
$\de^{old-gl}(x)=s_0x\otimes t_0x$ for $x\in\C_1$. This makes
$H_1^{old-gl}$ slightly change. It does not matter because there is no
influence on any potential applications. The difference appears in
a situation like that of
Figure~\ref{false}. With $C_0^{old-gl}(\C)=\Z\C_0\oplus \Z\C_0$,
$u+x-w-v$ is a old globular cycle. With
$C_0^{old-gl}(\C)=\Z(\C_0\otimes\C_0)$, this fake $1$-globular cycle is
killed.

\bth\label{relation-old_new}
We have the following commutative diagram of natural transformations
for $*\geq 0$ {\large
\[
\xymatrix{
&&{H_*^{gl}} \ar@{->}[rr]^{h^\pm}\ar@{->}[d]^{R^{gl}} &&  {H_*^\pm} \ar@{->}[d]^{R^{\pm}}\\
{H_*^{old-gl}}\ar@{->}[rru]\ar@/^40pt/[rrrru]^{(h^\pm)^{old}}\ar@{->}[rr]^{}\ar@{->}[rrd]^{}&& {\HR_*^{gl}} \ar@{->}[rr]^{h^\pm}
&&  {\HR_*^\pm} \\
&&{\HF_*^{gl}} \ar@{->}[u]_{\square^{gl}}\ar@{->}[rr]^{h^\pm} && {\HF_*^\pm}  \ar@{->}[u]_{\square^{\pm}}      }
\]}
where
\begin{itemize}
\item the map $H_*^{old-gl}\rightarrow H_*^{gl}$ is the canonical
map induced by $x\mapsto \square_n^{gl}(x)$ from $\C_n$ to
$\mathcal{N}_{n-1}^{gl}(\C)$
\item the map $H_*^{old-gl}\rightarrow \HF_*^{gl}$ is the canonical
map making all identifications like $A*_n B=A+B$ for any $n\geq 1$ and
any $p$-morphisms $A$ and $B$ with $p\geq n+1$
\item the map $\HF_*^{gl}\rightarrow \HF_*^{\pm}$  is the canonical
map making the supplemental identification $x=x*_0 y$ or $y=x *_0 y$
depending on the sign $\pm$
\item the map $\HF_*^{\pm}\rightarrow \HR_*^\pm$ is the canonical map
  induced by the folding operators $\square^\pm$ of \cite{Coin} (which
  is likely to be  an isomorphism for any strict globular
  $\omega$-category), and the map $\HF_*^{gl}\rightarrow \HR_*^{gl}$ is
  the canonical map induced by the folding operators $\square^{gl}$
  (which is also likely to be  an isomorphism for any strict globular
  $\omega$-category)
\item the maps $R^{gl,\pm}$ are the canonical maps from the globular
  or corner homology to the corresponding reduced homology (which are
  conjecturally an isomorphism for any free $\omega$-category
  generated by a semi-cubical set or a globular set).
\end{itemize}
\eth

\bpf This is due to the fact that for $n\geq 1$, the natural map
$(h^\pm_n)^{old}$ is induced by the set map $\square_n^-$ from $\C_n$
to $\omega Cat(I^n,\C)^-$ (\cite{Gau} Proposition 7.4). \epf

The difference between $H_0^{old-gl}$ and $H_0^{gl}$ is also not
important. The group $H_0^{old-gl}$ was indeed only introduced to
define the morphisms $h^-$ and $h^+$ in dimension $0$. But
$H_0^{old-gl}$ does not have any computer-scientific meaning and
is not involved in any potential applications.

\section{Globular homology and deformation of HDA}\label{deformation}

The following table summarizes how the globular nerve may be
understood and compared with the two corner nerves of $\C$.

{
\begin{center}
\begin{tabular}{|p{2.5cm}||p{2.5cm}|p{2.5cm}|p{2.5cm}|}
\hline
Geometric object & Formal theory & ``True''  theory & Simplicial cut\\\hline\hline
Branching & formal negative corner homology & negative corner homology & $\mathcal{N}^-(\C)$\\\hline
Merging & formal positive corner homology & positive corner homology &$\mathcal{N}^+(\C)$ \\\hline
Globe & formal globular homology & globular homology  & $\mathcal{N}^{gl}(\C)$\\\hline
\end{tabular}
\end{center}}

Intuitively, the globular nerve of $\C$ contains all
\textit{achronal cuts} in the middle of all globes, whereas the
negative and positive corner simplicial nerves contain all
\textit{achronal cuts} close to respectively the negative and
the positive corners of the automaton. The expression ``achronal'' is borrowed from
\cite{HDA2} and \cite{FGR}. In these papers, HDA are modeled by
local pospaces, and an achronal subspace $Y$ of a local pospace is a
topological subspace such that $x\leq y$ and $x,y\in Y$ imply $x=y$.
The remarkable point is that
the set of all achronal cuts of a given type can be enclosed
into a simplicial set.

This could mean that the whole geometry of the free $\omega$-category
$\C$ generated by a semi-cubical set (i.e. a HDA) would be contained in the
following diagram of augmented simplicial sets
\[\xymatrix{& {\mathcal{N}^{gl}(\C)}\ar@{->}[ld]_{h^-}\ar@{->}[rd]^{h^+}&\\
{\mathcal{N}^{-}(\C)}&&{\mathcal{N}^{+}(\C)}}\] and in its temporal
graph $tr^1\C$. This latter contains the information about the
temporal structure of the HDA.

A problem, already mentioned in \cite{ConcuToAlgTopo}, is the
question of the invariance of the globular homology of an
$\omega$-category up to a choice of a \textit{cubification}
\footnote{Some authors \cite{HDA} \cite{cyl} use the term
\textit{cubicalation} : this means decomposing a HDA in cubes.}
of the corresponding HDA.  There are two types of deformations :
the \textit{spatial deformations} or \textit{S-deformations} and
the \textit{temporal deformations} or \textit{T-deformations}.

The globular cut is invariant by S-deformation, that is by
deformations of $p$-morphisms with $p\geq 2$. This is simply due to
the fact that such a deformation corresponds in the globular cut to a
deformation of any simplex containing it as label. Therefore such a
deformation corresponds to a deformation up to homotopy, in the usual
sense, of the globular cut.

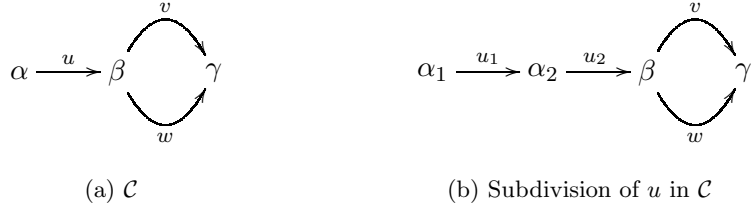
\begin{figure}
\begin{center}
\subfigure[$\C$]{\label{dilate1}
\xymatrix{\alpha\fr{u}& \beta \ar@/^20pt/[r]^v \ar@/_20pt/[r]_w&\gamma}}
\hspace{2cm}
\subfigure[Subdivision  of $u$ in $\C$]{
\label{dilate2}
\xymatrix{\alpha_1\fr{u_1}&\alpha_2 \fr{u_2}& \beta \ar@/^20pt/[r]^v \ar@/_20pt/[r]_w&\gamma }}
\end{center}
\caption{Subdivision of time}
\label{dilate}
\end{figure}

Unlike the corner homologies, the globular homology turns indeed to
depend on the subdivision of time. The reason is contained in
Figure~\ref{dilate}. The obvious $1$-functor from the left to the
right such that $u\mapsto u_1 *_0 u_2$ should leave the globular
homology invariant.  This is not the case because the first globular
homology is for the left member the free $\Z$-module generated by
$v-w$ and $u*_0v-u*_0w$, and for the right member the free $\Z$-module
generated by $v-w$ and $u_2*_0v-u_2*_0w$ and
$u_1*_0u_2*_0v-u_1*_0u_2*_0w$.  However in Figure~\ref{dilate}, one
can subdivide as many times as one wants for example $v$, and the
globular homology will not change.

One way to overcome this problem is exposed in the last sections
of \cite{ConcuToAlgTopo}, devoted to the description of a generic
way to produce T-invariants starting from the globular nerve.
Let us prove \cite{ConcuToAlgTopo} Claim 5.1 which enables to
introduce the bisimplicial set mentioned in that paper.

Let $\C$ be a non-contracting $\omega$-category. Using
Theorem~\ref{grading}, recall that for some $\omega$-functor $x$ from $\Delta^n$ to
$\P\C$, one calls $S(x)$ the unique element of the image of
$s_0\circ x$ and $T(x)$ the unique element of the image of
$t_0\circ x$. If $(\alpha,\beta)$ is a pair of $\mathcal{N}_{-1}^{gl}(\C)$,
set $S(\alpha,\beta)=\alpha$ and $T(\alpha,\beta)=\beta$.

\bp Let $\C$ be a non-contracting $\omega$-category. Let $x$
and $y$ be two $\omega$-functors from $\Delta^n$ to $\P\C$ with $n\geq
0$. Suppose that $T(x)=S(y)$.
Let $x*y$ be the map from the faces of
$\Delta^n$ to $\C$ defined by \[(x*y)((\sigma_0\dots
\sigma_r)):=x((\sigma_0\dots \sigma_r))*_0 y((\sigma_0\dots
\sigma_r)).\]
Then the following conditions are equivalent :
\begin{enumerate}
\item\label{c1} The image of $x*y$ is a subset of $\P\C$.
\item\label{c2} The set map $x*y$ yields an $\omega$-functor from $\Delta^n$ to $\P\C$
and $\de_i(x*y)=\de_i(x)*\de_i(y)$ for any $0\leq i\leq n$.
\end{enumerate}
On contrary, if for some $(\sigma_0\dots \sigma_r)\in \Delta^n$,
$(x*y)((\sigma_0\dots \sigma_r))$ is $0$-dimen\-sio\-nal, then $x*y$ is
the constant map $S(x)=T(y)$.
\ep

\bpf  We have to prove that Condition~\ref{c1} implies Condition~\ref{c2}. Let us
consider $P(n)$ : ``for any non-contracting $\omega$-category $\C$ and
any $\omega$-functor $x$ and $y$ from $\Delta^n$ to $\P\C$ such that
$T(x)=S(y)$ and such that the image of $x*y$ is a subset of $\P\C$,
then $x*y$ yields an $\omega$-functor from $\Delta^n$ to $\P\C$ and $\de_i(x*y)=\de_i(x)*\de_i(y)$
for any $0\leq i\leq n$.''

Property $P(0)$ is obvious. Suppose $P(n-1)$ proved for
$n\geq 1$. For any $0\leq i\leq n$, $\de_i(x)*\de_i(y)$ is a set map
from $\Delta^{n-1}$ to $\P\C$ satisfying the hypothesis of the
proposition, so by induction hypothesis, $\de_i(x)*\de_i(y)$ yields an
$\omega$-functor from $\Delta^{n-1}$ to $\P\C$. Let $z_i:=\de_i(x)*\de_i(y)$.
For $0\leq j<i\leq n$,
\begin{alignat*}{2}
\de_j (z_i)&= (\de_j \de_i(x))*(\de_j\de_i(y)) &&\ \hbox{ by induction hypothesis}\\
&= (\de_{i-1}\de_j (x))*(\de_{i-1}\de_j(y))&&\\
&= \de_{i-1} (\de_j(x)*\de_j(y)) &&\ \hbox{ by induction hypothesis}\\
&= \de_{i-1} z_j &&
\end{alignat*}
Therefore $(z_i)_{0\leq i\leq n}$ is an $(n-1)$-shell. So it provides
a unique $\omega$-functor
\[z:\Delta^n\backslash\{(01\dots n)\}\rightarrow \P\C\]
by Proposition~\ref{filling_simp}. It remains to check that
\[z\left(s_{n-1} R((01\dots n))\right)= s_{n}((x*y)((01\dots n)))\] and
\[z\left(t_{n-1} R((01\dots n))\right)= t_{n}((x*y)((01\dots n)))\] to
complete the proof. Let us check the first equality. One has
$$s_{n-1} R((01\dots n))=\Psi(X_1,\dots, X_s)$$ where $\Psi$ uses only
composition laws and where $X_1,\dots, X_s$ are faces of $\Delta^{n}$
of dimension at most $n-1$. Denote by $\Psi'$ the same function as $\Psi$ with
$*_i$ replaced by $*_{i+1}$. Then

{
\begin{alignat*}{2}
&z\left(s_{n-1} R((01\dots n))\right) &&\\
&= z \Psi(X_1,\dots, X_s) &&\\
&= \Psi'(z(X_1),\dots,z( X_s)) &&\ \hbox{ since $z$ $\omega$-functor}\\
&= \Psi'(x(X_1)*_0 y(X_1),\dots,x( X_s)*_0 y( X_s)) &&\ \hbox{ by definition of $z$}\\
&= \Psi'(x(X_1),\dots,x( X_s))*_0 \Psi'(y(X_1),\dots,y( X_s)) &&\ \hbox{ by interchange law}\\
&= \left(x \Psi(X_1,\dots, X_s)\right) *_0 \left(y \Psi(X_1,\dots, X_s)\right)&&\ \hbox{ since $x$ and $y$ $\omega$-functors}\\
&= \left(x s_{n-1} R((01\dots n))\right) *_0 \left(y s_{n-1} R((01\dots n))\right)&&\\
&= \left(s_{n} x  R((01\dots n))\right) *_0 \left(s_{n} y R((01\dots n))\right)&&\ \hbox{ since $x$ and $y$ $\omega$-functors}\\
&= s_{n}\left( x  R((01\dots n)) *_0  y R((01\dots n))\right)&&\ \hbox{ by interchange law}\\
&= s_{n}((x*y)((01\dots n))) &&
\end{alignat*}}

Now let us suppose that $(x*y)((\sigma_0\dots \sigma_r))$ is
$0$-dimensional in $\C$ for some $(\sigma_0\dots \sigma_r)$. Then
$$s_1 x((\sigma_0\dots \sigma_r)) *_0 s_1 y((\sigma_0\dots
\sigma_r))$$
is $0$-dimensional. Either $s_0(\sigma_0\dots
\sigma_r)=(n)$ (the initial state of $\Delta^n$) 
or there exists a $1$-morphism $U$ of $\Delta^n$ such
that $s_0U=(n)$ and $t_0U=s_0(\sigma_0\dots \sigma_r)$.  In the first
case, $x((n))*_0 y((n))$ is $0$-dimensional.  In the second case,
$$x(t_0U)*_0 y(t_0U)=t_1 x(U) *_0 t_1 y(U)= t_1\left( x(U) *_0
y(U)\right)$$ is $0$-dimensional. Then $x(U)*_0 y(U)$ is
$0$-dimensional as well as $$x((n))*_0 y((n))=s_1\left(x(U)*_0 y(U)\right).$$
For any face $(\tau_0\dots \tau_r)$ of $\Delta^n\backslash\{(n)\}$, there exists a
$1$-morphism $V$ from $((n))$ to $s_0(\tau_0\dots \tau_r)$ or
$t_0(\tau_0\dots \tau_r)$ : let us say $s_0(\tau_0\dots \tau_r)$.
Since $$s_1(x*y)(V)=(x*y)((n))$$ is $0$-dimensional, then
$(x*y)(V)$ is $0$-dimensional, as well as
$$t_1(x*y)(V)=(x*y)(s_0(\tau_0\dots \tau_r))=s_1(x*y)((\tau_0\dots \tau_r)).$$
Therefore $(x*y)((\tau_0\dots \tau_r))$ is $0$-dimensional.
\epf

In the sequel, we set $(\alpha,\beta)*(\beta,\gamma)=(\alpha,\gamma)$,
$S(\alpha,\beta)=\alpha$ and $T(\alpha,\beta)=\beta$.
If $x$ is an $\omega$-functor
from $\Delta^n$ to $\P\C$, and if $y$ is the constant map $T(x)$ (resp. $S(x)$)
from $\Delta^n$ to $\C_0$, then set $x*y:=x$ (resp. $y*x:=x$).

\bth Suppose that $\C$ is an object of $\omega Cat_1$. Then for
$n\geq 0$, the operations $S$, $T$ and $*$ allow to define a
small category $\underline{\mathcal{N}^{gl}_n(\C)}$ whose
morphisms are the elements of $\mathcal{N}^{gl}_n(\C)\cup
\{\hbox{constant maps }\Delta^n\rightarrow \C_0\}$ and whose
objects are the $0$-morphisms of $\C$. If
$\underline{\mathcal{N}^{gl}_{-1}(\C)}$ is the small category
whose morphisms are the elements of $\C_0\p\C_0$ and whose
objects are the elements of $\C_0$ with the operations $S$, $T$
and $*$ above defined, then one obtains (by defining the face
maps $\de_i$ and degeneracy maps $\epsilon_i$ in an obvious way
on $\{\hbox{constant maps }\Delta^n\rightarrow \C_0\}$) 
an augmented simplicial object $\underline{\mathcal{N}^{gl}_*}$
in the category of small categories. \eth

\bpf
Equalities     $S(x)=\de_i S(x)$, $S(x)=\epsilon_i S(x)$,
$T(x)=\de_i T(x)$, $T(x)=\epsilon_i T(x)$ are consequences of
Proposition~\ref{grading}.  Equality $\de_i(x*y)=\de_i x * \de_i y$ is
proved right above.  The verification of
$\epsilon_i(x*y)=\epsilon_i x * \epsilon_i y$ is straightforward.
\epf

By composing by the classifying space functor of small categories
(cf. for example \cite{classifiant} for further details), one obtains
a bisimplicial set which seems to be well-behaved with respect to
subdivision of time.  Indeed the first total homology groups
associated to both $\omega$-categories of Figure~\ref{dilate} are
equal to $\Z$. Further explanations will be given in future papers.

To conclude, let us point out that in reasonable cases, i.e. when the
$p$-morphisms (with $p\geq 2$) of a non-contracting $\omega$-category
$\C$ are invertible with respect to the composition laws $*_i$ of $\C$
for $i\geq 1$, then $\P\C$ becomes a globular $\omega$-groupoid in the
sense of Brown-Higgins. And therefore in such a case, it is well-known
that the globular nerve of $\C$ satisfies the Kan property (see
\cite{oriental} or a generalization in \cite{MR89m:18011}).  However,
this is not true in general for both corner nerves. To understand this
fact, consider the $2$-source of $R(000)$ in Figure~\ref{I3} and
remove $R(0+0)$.  Consider both inclusion $\omega$-functors from $I^2$
to respectively $R(-00)$ and $R(00-)$. Then the Kan condition fails
because one cannot make the sum of $R(-00)$ and $R(00-)$ since
$R(0+0)$ is removed.


{\ \\
Institut de Recherche Math\'ematique Avanc\'ee\\
ULP et CNRS\\
7 rue Ren\'e Descartes\\
67084 Strasbourg Cedex\\
France\\
gaucher@irma.u-strasbg.fr}

\end{document}